\pdfoutput=1
\documentclass[a4wide]{article}

\usepackage{graphicx}

\usepackage[utf8]{inputenc}
\usepackage[english]{babel}
\graphicspath{ {images/} }
\usepackage{algorithm2e}
\usepackage{color}
\usepackage {amsmath,amsfonts,amssymb}
\usepackage{subfig}
\usepackage{color}

\usepackage{texnames}
\usepackage{url}
\usepackage{graphicx}
\definecolor{hotpink}{rgb}{0.9,0,0.5}
\usepackage[colorlinks,urlcolor=blue,citecolor=hotpink,linkcolor=blue]{hyperref}
\def\ycite[#1#2#3#4#5]#6{\cite[$\mit{#1#2#3#4}$#5]{#6}}

\newcommand{\R}{\mathbb{R}}
\newcommand{\F}{\mathbb{F}}
\newcommand{\fl}{\mathrm{fl}}

\renewcommand{\leq}{\leqslant}
\renewcommand{\geq}{\geqslant}

\begin{document}

\title{Reducing and monitoring round-off error propagation for symplectic implicit Runge-Kutta schemes
}

%\author[1]{Mikel Anto\~nana}
%\author[1]{Joseba Makazaga}
%\author[1]{Ander Murua}
%\affil[1]{
\author{Mikel Anto\~nana, Joseba Makazaga,  Ander Murua \\
KZAA saila, Informatika Fakultatea, UPV/EHU\\ Donostia / San Sebasti\'an}

\date{  }

\maketitle

\begin{abstract}
We propose an implementation of symplectic implicit Runge-Kutta schemes 
for highly accurate numerical integration of non-stiff Hamiltonian systems based on fixed point iteration.
Provided that the computations are done in a given floating point arithmetic, the precision of the results is limited by round-off error propagation.   We claim that our implementation with fixed point iteration is near-optimal with respect to round-off error propagation under the assumption that  the function that evaluates the right-hand side of the differential equations is implemented with machine numbers (of the prescribed floating point arithmetic) as input and output. In addition, we present a simple procedure to estimate the round-off error propagation by means of a slightly less precise second numerical integration. 
Some numerical experiments are reported to illustrate
the round-off error propagation properties of the proposed implementation.

\end{abstract}

\section{Introduction}
\label{sec:intro}

When numerically integrating an autonomous Hamiltonian system, one typically monitors the error in the preservation of the Hamiltonian function to check the precision of the numerical solution. However, severe loss of precision can actually occur 
for sufficiently long integration intervals, while displaying a good preservation of the value of Hamiltonian function. For high precision numerical integrations, where round-off errors may dominate truncation errors,  it is highly desirable both reducing and monitoring the propagation of round-off errors.

We propose an implementation of symplectic implicit Runge-Kutta schemes (such as RK collocation methods with Gaussian nodes) that takes special care in reducing the propagation of round-off errors. Our implementation is intended to be applied to non-stiff problems, which motivates us to solve the implicit equations by fixed-point iteration (see for instance~\cite{JMSanz-Serna1994}\cite{Hairer2006} for numerical tests comparing the efficiency of implementations based on fixed point iterations  and simplified Newton).

We work under the assumption that the (user defined) function that evaluates the right-hand side of the differential equations is implemented in such a way that input and output arguments are machine numbers of some prescribed floating point arithmetic.
Our actual implementation includes the option of computing, in addition to the numerical solution, an estimation of the propagated round-off error.

The starting point of our implementation is the work of Hairer et al.~\cite{Hairer2008}. There, the authors observe that a standard fixed-point implementation of symplectic implicit RK (applied with compensated summation~\cite{Higham2002}) exhibits an unexpected systematic error in energy due to round-off errors, not observed in explicit symplectic methods. They make the following observations that allow them to understand  that unfavorable error behavior:  (a) The implicit Runge-Kutta method whose coefficients $\tilde b_i,\tilde a_{i j}$ are the floating-point representation of the coefficient $b_i,a_{i,j}$ of a symplectic Runge-Kutta method is not symplectic; (b)  The error due to the application at each step of a fixed point iteration with standard stopping criterion (depending on a prescribed tolerance of the iteration error) tends to have a systematic character.  Motivated by these observations, they modify the standard implementation of fixed point iteration which allows them to reduce the effect of round-off errors. No systematic error in energy is observed in the numerical experiments reported in~\cite{Hairer2008}. However, we observe in some numerical experiments that the stopping criterion for the fixed point iteration that they propose fails to work properly in some cases. In addition, we claim that their implementation is still not optimal with respect to round-off error propagation.

In Section~3, we propose alternative modifications of the standard fixed point implementation of symplectic implicit Runge-Kutta methods, which compare very favorably with that proposed in~\cite{Hairer2008}. 

We first define a reference implementation with fixed point iteration where all the arithmetic operations  other than the evaluation of the right-hand side of the system of differential equations are performed in exact arithmetic, and as many iterations as needed are performed in each step. Such an implementation, that we call  FPIEA (Fixed Point iteration with Exact Arithmetic) implementation, is based on the following two modifications to the standard implementation with fixed point iterations: 
(i) From one hand, we reformulate each symplectic implicit Runge-Kutta method in such a way that its coefficients can be approximated by machine numbers while still keeping its symplectic character exactly (Subsection~3.1).  (ii) On the other hand, 
we propose a modification of the stopping criterion introduced in~\cite{Hairer2008} that is more robust and is independent of the chosen norm (Subsection~3.2)

The implementation we present here is based on the FPIEA implementation, with most multiplications and additions performed (for efficiency reasons) in the prescribed floating point arithmetic, but some of the operations performed with special care in order to reduce the effect of round-off errors. In particular, this includes a somewhat non-standard application of  Kahan's "compensated summation" algorithm~\cite{Kahan1965}\cite{Higham2002}\cite{Muller2009}, described in detail in Subsection~3.3.

Finally, in Subsection~3.4,  we present a simple procedure to estimate the round-off error propagation as the difference of the actual numerical solution, and a slightly less precise second numerical solution.
These two numerical solutions can be computed either in parallel, or sequentially with a lower computational cost than two integrations executed in completely independent way.

In Section~4, we show some numerical experiments to asses the performance of our final implementation. Some concluding remarks are presented in Section~5.

\section{Preliminaries}

\subsection{Numerical integration of ODEs by symplectic IRK schemes}

We are mainly interested in the application of symplectic implicit Runge-Kutta (IRK) methods for the numerical integration of Hamiltonian systems of the form 
\begin{equation}
\label{eq:Hamsyst}
\frac{d}{dt}q^j = \frac{\partial H(p,q)}{\partial p^j},
\quad 
\frac{d}{dt}p^j = -\frac{\partial H(p,q)}{\partial q^j}, \quad j=1,\ldots,d,
\end{equation}
where $H:\R^{2d} \to \R$. Recall that the Hamiltonian function $H(q,p)$ is a conserved quantity of the system.

More generally, we consider initial value problems of systems of autonomous ODEs of the form
\begin{equation}
\label{eq:ivp}
\frac{d}{dt}y=f(y),\quad  y(t_0)=y_0,
\end{equation}
where $f: \R^D\to \R^D$ is a sufficiently smooth map and $y_0 \in \R^D$. In the case of the Hamiltonian system (\ref{eq:Hamsyst}), $D=2d$, $y=(q^1,\ldots,q^d,p^1,\ldots,p^d)$.

For the system of differential equations (\ref{eq:ivp}), an s-stage implicit Runge-Kutta method is determined by an integer $s$  and the real coefficients $a_{ij}$ ($1 \leq i, j \leq s$), $b_{i}$ ($1\leq i \leq s$). The approximations $y_{n} \approx y(t_n)$ to the solution $y(t)$ of (\ref{eq:ivp}) at $t=t_{n}=t_0 + n h$ for $n=1,2,3,\ldots$ are computed as
\begin{equation}
\label{eq:yz}
y_{n+1}=y_n+h\sum^s_{i=1} b_i \, f(Y_{n,i}),
\end{equation}   
where the {\em stage vectors} $Y_{n,i}$ are implicitly defined at each step by
\begin{equation}
\label{eq:Y}
Y_{n,i} =y_n+ h \sum^s_{j=1}{a_{ij}\,f(Y_{n,j})}, \quad  i=1 ,\ldots, s.
\end{equation}
An IRK scheme is symplectic if an only if~\cite{JMSanz-Serna1994}
\begin{equation} \label{eq:sympl_cond_1}
b_{i}a_{ij}+b_{j}a_{ji}-b_{i}b_{j}=0, \ \ 1 \leq i,j \leq s.
\end{equation}
In that case, the IRK scheme  conserves exactly all quadratic first integrals of the original system (\ref{eq:ivp}), and if the system is Hamiltonian,  under certain assumptions~\cite{Hairer2006}, it approximately conserves the value of the Hamiltonian function $H(y)$ over long time intervals.

\subsection{Floating point version of an IRK integrator}

Let $\F \subset \R$ be the set of machine numbers of a prescribed floating point system. Let $\fl:\R \longrightarrow \F$ be a map that sends each real number $x$ to a nearest machine number $\fl(x) \in \F$.

We assume that instead of the original map $f: \R^D \to \R^D$, we have a computational substitute 
\begin{equation}
\label{eq:tildef}
\tilde f: \F^D\to \F^D.
\end{equation}
Ideally, for each $x \in \F^D$, $\tilde f(x):= \fl(f(x))$. In practice, the intermediate computations to evaluate $\tilde f$ are typically made using the floating point arithmetic corresponding to $\F$, which will result in some error $||\tilde f(x) - \fl(f(x))||$ caused by the accumulated effect of several round-off errors. 

We aim at efficiently implementing a given symplectic IRK scheme under the assumption that $f: \R^D \to \R^D$ is replaced by  (\ref{eq:tildef}).  Hence, the effect of round-off errors will be present even in the best possible ideal implementation where exact arithmetic were used for all the computations except for the evaluations of the map (\ref{eq:tildef}).  Our goal is to implement the IRK scheme working at the prescribed floating point arithmetic, in such a way that the effect of round-off errors is similar in nature and relatively close in magnitude to that of such ideal implementation.

\subsection{Kahan's compensated summation algorithm}

Obtaining the numerical approximation  $y_{n} \approx y(t_{n})$, ($n=1,2,\ldots$) to the solution $y(t)$ of the initial value problem (\ref{eq:ivp}) defined by (\ref{eq:yz})--(\ref{eq:Y}) requires computing the sums
\begin{equation}
\label{eq:sumy_n}
y_{n+1} = y_{n} + x_n, \quad n=0,1,2,\ldots, 
\end{equation}
where 
\begin{equation*}
x_n = h\sum^s_{i=1} b_i \, f(Y_{n,i}).
\end{equation*}

For an actual implementation that only uses a floating point arithmetic with machine numbers in $\F$, special care must be taken with the additions (\ref{eq:sumy_n}). It is worth mentioning that for sufficiently small step-length $h$,   the components of $x_n$ are smaller in size than those of $y_{n}$ (provided that the components of the solution $y(t)$ of (\ref{eq:ivp}) remain away from zero).   The naive recursive algorithm 
$\hat y_{n+1} :=\fl(\hat y_{n} + \fl(x_n))$,  ($n=0,1,2,3\ldots$),
typically suffers, for large $n$,  a significant loss of precision due to round-off errors.  It is well known that such a round-off error accumulation can be greatly reduced with the use of Kahan's compensated summation algorithm~\cite{Kahan1965} (see also~\cite{Higham2002},~\cite{Muller2009}). 
  Given a sequence $\{y_0,x_0,x_1,\ldots,x_n,\ldots\} \subset \F$ of machine numbers, Kahan's algorithm is aimed to compute the sums $y_n = y_0 + \sum_{\ell=0}^{n-1} x_{\ell}$, ($n\geq 1$,) using a prescribed floating point arithmetic, more precisely than with the naive recursive algorithm. In Kahan's algorithm, machine numbers $\tilde y_n$ representing the sums $y_n$ are computed along with additional machine numbers $e_n$ intended to capture the error $y_n-\tilde y_n$. The actual algorithm reads as follows: 
 \begin{algorithm}[H]
  \BlankLine
   $\tilde y_0= y_0; \ e_0=0$\;
   \BlankLine
   \For{$l\leftarrow 0$ \KwTo $n$}
   {
    \BlankLine
     $X_l = \fl(x_l + e_{l})$\;
     $\tilde y_{l+1} = \fl(\tilde y_{l} + X_{l})$\;
     $\hat X_{l} = \fl(\tilde y_{l+1} - \tilde y_{l})$\;
     $e_{l+1} =  \fl(X_{l} - \hat X_{l})$\;
    \BlankLine
   }
  \caption{Kahan’s compensated summation.}
  \label{alg:Kahan'sCS}
 \end{algorithm}
 The sums $\tilde y_l + e_l$ (which in general do not belong to $\F$) are more precise approximations of the exact sums $y_l$ than $\tilde y_l \in \F$. In this sense, if $y_0 \not \in \F$, the algorithm (\ref{alg:Kahan'sCS}) should be initialized as $\tilde y_0:= \fl(y_0)$ and $e_0:=\fl(y_0-\tilde y_0)$ (rather than $e_0=0$). 

Of course, algorithm (\ref{alg:Kahan'sCS}) also makes sense for $D$-vectors of machine numbers $\tilde y_0, e_0,x_0,x_1,\ldots,x_n \in \F^D$.
In this setting, algorithm (\ref{alg:Kahan'sCS}) can be interpreted as a family of maps parametrized by $n$ and $D$,
\begin{equation*}
S_{n,D}:\F^{(n+3)D} \to \F^{2D},
\end{equation*}
that given the arguments $\tilde y_0,e_0,x_0,x_1,\ldots,x_n \in \F^D$,  returns $\tilde y_{n+1}, e_{n+1} \in \F^D$ such that $\tilde y_{n+1} + e_{n+1}$ is intended to represent the sum $(\tilde y_0 + e_0) + x_0 +x_1 + \cdots + x_n$ (with some small error).

\section{Proposed implementation of symplectic IRK schemes}
\label{sec:2}

\subsection{Symplectic schemes with machine number coefficients}
\label{ss:3.1}

If the coefficients  $b_i,a_{i j}$ determining a symplectic IRK are replaced by machine numbers $\tilde b_i,\tilde a_{i j} \in \F$ that approximate them (say, $\tilde b_i := \fl(b_i)$, $\tilde a_{i j} := \fl(a_{i j})$), 
then the resulting IRK scheme typically fails to satisfy the symplecticity conditions (\ref{eq:sympl_cond_1}). This results in a method that does not conserve quadratic first integrals and exhibits a linear drift in the value of the Hamiltonian function~\cite{Hairer2008}.

Motivated by that, we recast the definition of a step of the IRK method  as follows:
\begin{align}
\label{eq:YL}
Y_{n,i}  &=y_n+ \sum^s_{j=1}{\mu_{ij}\,L_{n,j}, \quad  L_{n,i} = h b_i f(Y_{n,i})}, \quad  i=1 ,\ldots, s, \\
\label{eq:y}
y_{n+1} &=y_n+\sum^s_{i=1} L_{n,i},
\end{align}   
where
\begin{equation*} 
\mu_{ij}=a_{ij}/b_j,  \quad 1 \leq i,j \leq s.
\end{equation*}
Condition (\ref{eq:sympl_cond_1}) now becomes,
\begin{equation*} 
\mu_{ij}+\mu_{ji}-1=0, \quad 1 \leq i,j \leq s.
\end{equation*}

The main advantage of the proposed formulation over the standard one is that the absence of multiplications in the alternative symplecticity condition  makes possible (see Appendix~A for the particular case of the 12th order Gauss collocation IRK method) to find machine number approximations $\tilde{\mu}_{i j}$ of $\mu_{i j}=a_{i j}/b_j$ satisfying exactly the symplecticity condition
\begin{equation} \label{eq:sympl_cond_2}
\tilde{\mu}_{ij}+\tilde{\mu}_{ji}-1=0, \quad 1 \leq i,j \leq s.
\end{equation}

\subsection{Iterative solution of the nonlinear Runge-Kutta equations} 
\label{ss:3.2}

The fixed point iteration can be used to approximately compute the
solution of the implicit equations (\ref{eq:YL}) as follows: For $k=1,2,\ldots$ obtain the approximations $Y_{n,i}^{[k]}$, $L_{n,i}^{[k]}$ of $Y_{n,i}$, $L_{n,i}$ ($i=1,\ldots,s$) as 
\begin{equation}
\label{eq:fixed_point_iteration}
L_{n,i}^{[k]} = h b_i\,  f(Y_{n,i}^{[k-1]}), \quad Y_{n,i}^{[k]}  =y_n+ \sum^s_{j=1}\,  \mu_{ij}\, L_{n,j}^{[k]} \quad  i=1 ,\ldots, s.
\end{equation}
The iteration may be initialized simply with $Y_{i}^{[0]}=y_n$, or by some other procedure that uses the stage values of the previous steps~\cite{Hairer2006}. If the step-length $h$ is sufficiently small, these iterations converge to a fixed point that is solution of the algebraic equations (\ref{eq:YL}). 

The situation is different for an actual computational version of these iterations, where $f$ is replaced in (\ref{eq:fixed_point_iteration})  by its computational substitute (\ref{eq:tildef}). The $k$th iteration then reads as follows: For $ i=1 ,\ldots, s$,
\begin{equation}
\label{eq:fixed_point_iteration_ideal}
f_{n,i}^{[k]} =  \tilde f({\fl(Y_{n,i}^{[k-1]})}), \quad L_{n,i}^{[k]} = h b_i\, f_{n,i}^{[k]}, \quad Y_{n,i}^{[k]}  = 
{ y_n+ \sum^s_{j=1}\,  \mu_{ij}\, L_{n,j}^{[k]} }.
\end{equation}
In this case,  either a fixed point of (\ref{eq:fixed_point_iteration_ideal}) is reached in a finite number of iterations, or the iteration fails (mathematically speaking) to converge. In the former case, however (provided that $h$ is  small enough for the original iteration (\ref{eq:fixed_point_iteration}) to converge),  after a finite number of iterations,  a computationally acceptable approximation to the fixed point of (\ref{eq:fixed_point_iteration}) is typically achieved, and the successive iterates remain close to it. According to our experience and the numerical experiments reported in~\cite{Hairer2008},   a computational fixed point is reached for most steps in a numerical integration with sufficiently small step-length $h$.

In standard implementations of implicit Runge-Kutta methods, one considers 
\begin{equation*}
\Delta^{[k]} = (Y_{n,1}^{[k]}-Y_{n,1}^{[k-1]}, \ldots,Y_{n,s}^{[k]}-Y_{n,s}^{[k-1]}) \in \F^{s d}, 
\end{equation*}
(for notational simplicity, we avoid reflecting the dependence of $n$ on $\Delta^{[k]}$), and stops the iteration provided that $||\Delta^{[k]}|| \leq \mathrm{tol}$, with a prescribed vector norm $||\cdot||$ and iteration error tolerance $\mathrm{tol}$. If the chosen value of $\mathrm{tol}$ is too small, then the iteration may never end when the computational sequence does not arrive to a fixed point in a finite number of steps. 
If $\mathrm{tol}$ is not small enough, the iteration will stop too early, which will result in an iteration error of larger magnitude than round-off errors. Furthermore, as observed in~\cite{Hairer2008}, such iteration errors tend to accumulate in a systematic way.

The remedy proposed in~\cite{Hairer2008} is to stop the iteration either if $\Delta^{[k]}=0$ (that is, if a fixed point is reached) or if $||\Delta^{[k]}|| \geq || \Delta^{[k-1]}||$. %
 The underlying idea is that (provided that $h$ is  small enough for the original iteration (\ref{eq:fixed_point_iteration}) to converge), typically $||\Delta^{[k]}|| < || \Delta^{[k-1]}||$ whenever the iteration error is substantially larger than round-off errors,  and thus $||\Delta^{[k]}|| \geq || \Delta^{[k-1]}||$ may indicate that round-off errors are already significant. 
 
\begin{figure}[h!]
\centering
\begin{tabular}{c c}
\subfloat[OSS: Hairer's stopping criterion]
{\includegraphics[width=.4\textwidth]{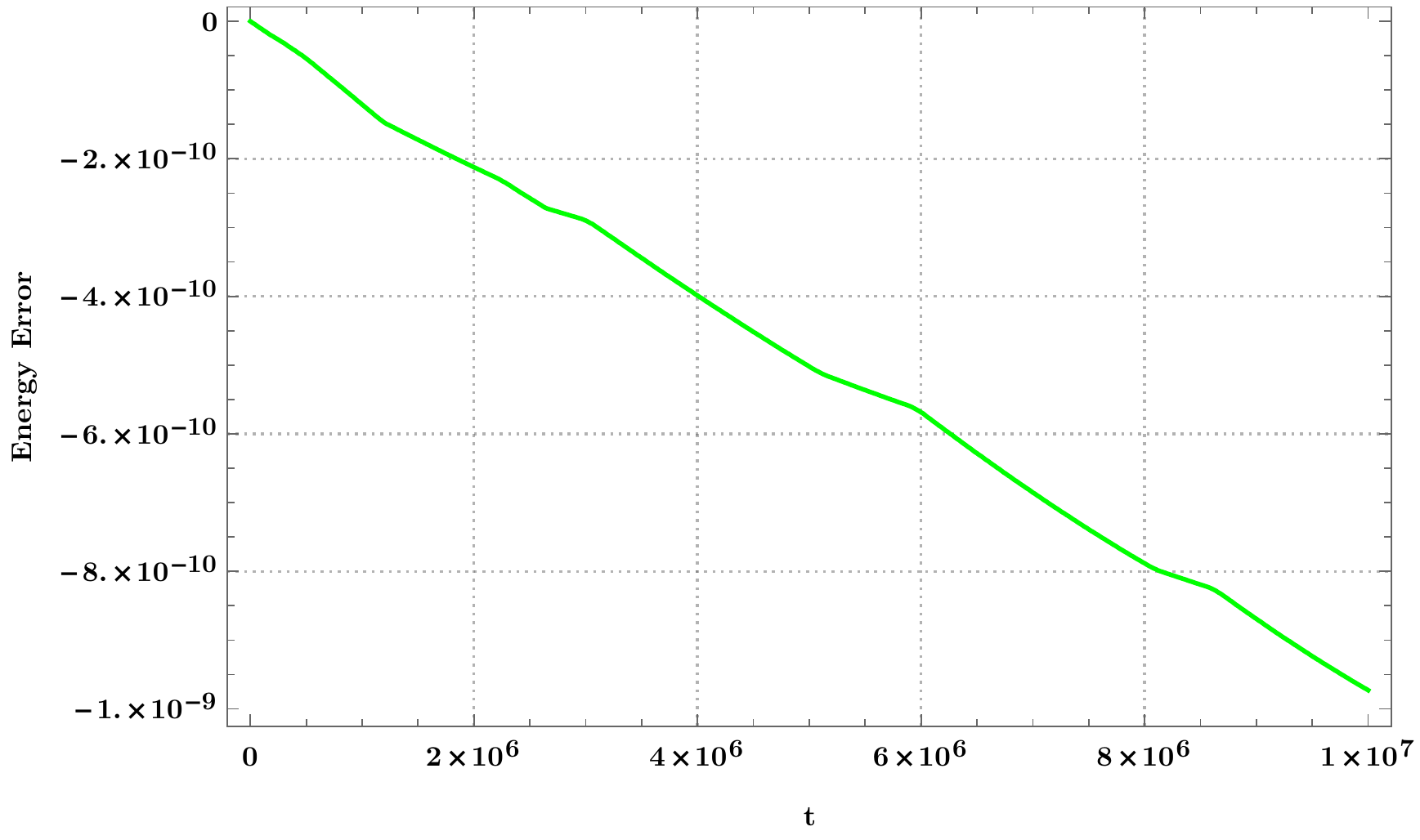}}
&
\subfloat[OSS: New stopping criterion]
{\includegraphics[width=.4\textwidth]{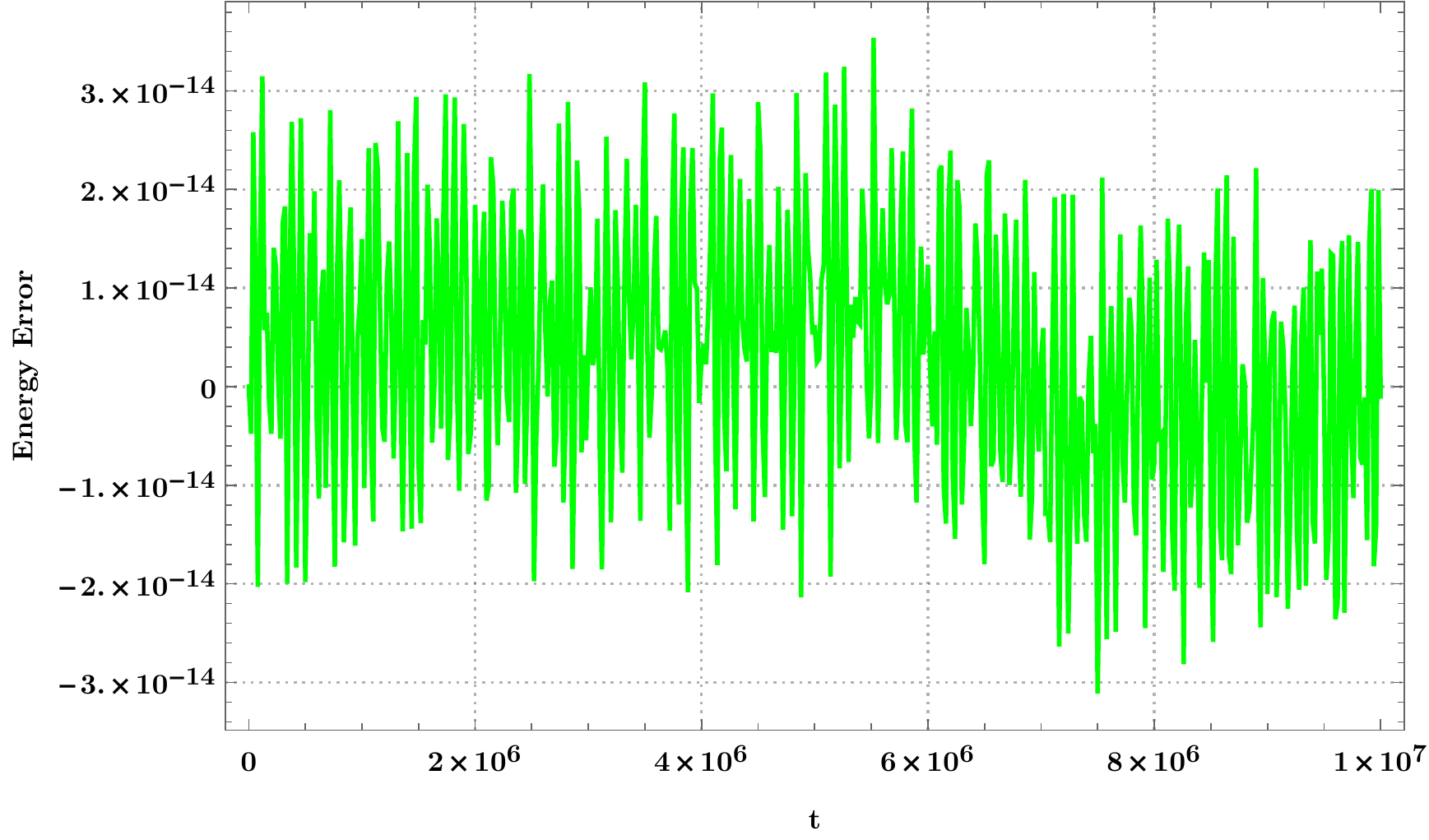}}
\end{tabular}
\caption{\small Evolution of relative error in energy for the outer solar system problem (OSS) with the original unperturbed initial values in~\cite{Hairer2008} and doubled step-size ($h=1000/3$ days).  (a) Hairer's stopping criterion. (b) New stopping criterion.}
\label{fig:plot0}
\end{figure}

We have observed that Hairer's strategy works well in general, but in some cases it stops the iteration too early. Indeed, it works fine for the initial value problem on a simplified model of the outer solar system (OSS) reported in~\cite{Hairer2008} with a step-size $h$ of $500/3$ days, but it goes wrong with $h=1000/3$.
Actually,  we have run   \href{http://www.unige.ch/~hairer/preprints/code.tar}{Hairer's fortran code} and observed that the computed numerical solution exhibits an error in energy that is considerably larger than round-off errors.  The evolution of relative error in energy is displayed on the left of Figure~\ref{fig:plot0}, which shows a linear growth pattern. We have checked that, for instance, at the first step,  
\begin{equation*}
\|\Delta^{[1]}\| > \|\Delta^{[2]}\|> \cdots >\|\Delta^{[12]}\| = 3.91\times 10^{-14} \le \|\Delta^{[13]}\| = 4.35\times 10^{-14}
\end{equation*}
 which causes the iteration to stop at the $13$th iteration, which happens to be too early,  since  subsequently, $\|\Delta^{[13]}\|> \|\Delta^{[14]}\| > \|\Delta^{[15]}\|> \|\Delta^{[16]}\|=0$.
% shows using his implementation. In the next section, we will repeat the same experiment using our implementation  (Fig.\ref{fig:plot1}) and this will keep a minimal growth of energy error.   

Motivated by that we propose an alternative more stringent stopping criterion: Denote as $\Delta_j^{[k]}$ the $j$th component ($1\leq j \leq s D$) of  $\Delta^{[k]} \in \R^{s D}$.  The fixed point iteration (\ref{eq:fixed_point_iteration_ideal}) is performed for $k=1,2,\ldots$ until either $ \Delta^{[k]} =0$ or the following condition is fulfilled for two consecutive iterations:
\begin{equation}
\label{eq:stopping}
\forall j \in \{1,\ldots,s D\},  \quad
\min \left(\{|\Delta_j^{[1]}|,\cdots ,|\Delta_j^{[k-1]}|\} \ /\{0\} \right) \leq |\Delta_j^{[k]}|.
\end{equation}
If $K$ is the first positive integer such that (\ref{eq:stopping}) does hold for both $k=K-1$ and $k=K$, then we compute the approximation $y_{n+1} \approx y(t_{n+1})$ as 
 \begin{equation*}
y_{n+1} = y_{n} + \sum_{i=1}^{s} L_{n,i}^{[K]}.
\end{equation*}
The iteration typically stops with $\Delta_j^{[K]}=0$ for all $j$. However, when the iteration stops with $\Delta_j^{[K]}\neq 0$ for some $j$, that is, when no computational fixed point is achieved,  we still have to decide if this has been due to the effect of small round-off errors, or because the step-size is not small enough for the convergence of the iteration (\ref{eq:fixed_point_iteration}). Here is the only point where our implementation depends on a norm-based standard criterion (with rather loose absolute and relative error tolerances) to decide if $\Delta^{[K]} \in \F^{s D}$ is small enough.

We have repeated the experiment of OSS with $h=1000/3$ by replacing Hairer's stopping criterion by our new one. The evolution of the resulting energy errors are displayed on the right of Figure~\ref{fig:plot0}.

\subsection{Machine precision implementation of new algorithm}

Subsections~\ref{ss:3.1} and~\ref{ss:3.2} completely determine the FPIEA (Fixed Point Iteration with Exact Arithmetic) implementation referred to in the Introduction. We next describe in detail our machine precision implementation of the  algorithm described (for exact arithmetic) in previous subsection.
 
 Consider appropriate approximations $\tilde b_i \in \F$ of $b_i$ ($i=1,\ldots,s$), and let 
 $\tilde{\mu}_{i j} \in \F$ ($i,j=1,\ldots,s$) be approximations of $\mu_{i j}$ satisfying exactly the symplecticity condition (\ref{eq:sympl_cond_2}).

Given $y_n \in \R^D$, we consider $\tilde y_0 = \fl(y_0)$ and $e_0=\fl(y_n-\tilde y_n)$. 
For each $n=0,1,2,\ldots$, we initialize $Y_{n,i}^{[0]}=y_n$, 
and successively compute for $k=1,2,\ldots$
\begin{equation}
\label{eq:fixed_point_iteration_comp}
\begin{split}
    f_{n,i}^{[k]} &= \tilde f(Y_{n,i}^{[k-1]}), \quad L_{n,i}^{[k]} = \fl(h\,  \tilde{b}_i\,f_{n,i}^{[k]}), \\
Z_{n,i}^{[k]} &= {e_n +} \sum^s_{j=1}\, \tilde{\mu}_{ij}\, L_{n,j}^{[k]}, \quad Y_{n,i}^{[k]}  = \fl\Big( \tilde y_n+ Z_{n,i}^{[k]}\Big)
\end{split}
\end{equation}
until the iteration is stopped at $k=K$ according to the criteria described in Subsection~\ref{ss:3.2}. Hence, $K$ is the highest index $k$ such that $f_{n,i}^{[k]}$ has been computed. 

In our actual implementation, one can optionally initialize $Y_{n,i}^{[0]}$ by interpolating from the stage values of previous step, which improves the efficiency of the algorithm. Nevertheless, in all the numerical  results reported Section~\ref{s:ne} below, the simpler initialization 
 $Y_{n,i}^{[0]}=y_n$ is employed.
 
In (\ref{eq:fixed_point_iteration_comp}),  we evaluate each $Z_{n,i}^{[k]} \in \F^D$  as 
\begin{equation*}
Z_{n,i}^{[k]} = (\cdots ((e_n + \tilde{\mu}_{i,1} L_{n,1}^{[k]}) + \tilde{\mu}_{i,2} L_{n,2}^{[k]} )+\cdots+\tilde{\mu}_{i,s-1} L_{n,s-1}^{[k]})+\tilde{\mu}_{i,s}L_{n,s}^{[k]}
\end{equation*}
where each multiplication and addition is performed in the prescribed floating point arithmetic. 

We then compute $\tilde y_{n+1}, e_{n+1} \in \F^D$ such that $\tilde y_{n+1} + e_{n+1} \approx y(t_{n+1})$
as follows: 
\begin{itemize}
\item compute for $i=1,\ldots,s$ the vectors   
\begin{equation}
\label{eq:Eni}
E_{n,i} = h\,  \tilde b_i\,f_{n,i}^{[K]}-L_{n,i}^{[K]}.
\end{equation}
\begin{equation*}
{\delta_{n}=e_{n}+\sum_{i=1}^{s} E_{n,i}}.
\end{equation*}
\item finally, compute 
  \begin{equation}
  \label{eq:y3}
    (\tilde y_{n+1}, e_{n+1})=S_{s,D}(\tilde y_{n},\delta_{n},L_{n,1}^{[K]}, \ldots, L_{n,s}^{[K]}).
  \end{equation}
\end{itemize}
If the FMA (fused-multiply-add) instruction is available, it should be used to compute (\ref{eq:Eni}) (with precomputed coefficients $h \tilde b_i$). The order in which the terms defining $Z_n$ and $\delta_n$ are actually computed in the floating point arithmetic is not relevant, as the corresponding round-off errors of the small corrections {$Z_n$} and $\delta_n$ will have a very marginal effect in the computation of $\tilde y_{n+1}$ and $e_{n+1}$.

\subsection{Round-off error estimation.}
\label{ss:estimation}

We estimate the round-off error propagation of our numerical solution $\tilde{y}_{n}+ e_n \approx y(t_n)$ ($n=1,2,\ldots$) by computing its difference with a slightly less precise secondary numerical solution $\hat y_{n}+\hat e_n \approx y(t_n)$ obtained with a modified version of the machine precision algorithm described in previous section. In this modified version of the algorithm,
the components of each $L_{n,i}^{[K]}$ in (\ref{eq:y3}) are rounded to a machine number with a shorter mantissa. We next give some more details. 

 Let $p$ be the number of binary digits of our floating point arithmetic.
Given an integer $r\geq 0$ and a machine number $x$, we define
$\fl_{p-r}(x):= \fl(2^r x + x) - 2^r x$. This is essentially equivalent to rounding $x$ to a floating point number with $p-r$ significant binary digits. 

We determine the algorithm for the secondary integration by fixing a positive integer $r<p$ and modifying (\ref{eq:y3}) in the implementation of the algorithm described in previous subsection as follows:
\begin{equation*}
  (\hat y_{n+1}, \hat e_{n+1})=S_{s,D}(\hat y_{n},\delta_{n},\fl_{p-r}(L_{n,1}^{[K]}), \ldots, \fl_{p-r}(L_{n,s}^{[K]})).
\end{equation*}

The proposed round-off error estimation can thus be obtained as the difference of the primary numerical solution and the secondary numerical solution obtained with a relatively small value of $r$ (say, $r=3$). These two numerical solutions can be computed in parallel in a completely independent way. 

In addition, we have implemented a sequential version
 with lower CPU requirements than two integrations executed sequentially in completely independent way. The key to do that is the following: At each step, the stage values $Y_{n,i}$ ($i=1,\ldots,s$) of both primary and secondary integration will typically be very close to each other (as far as the estimated round-off error does not grow too much). Thus, the number of iterations of each step of the secondary integration can be reduced by using the final stage values $Y_{n,i}$ ($i=1,\ldots,s$) of the primary integration as initial values $Y_{n,i}^{[0]}$ of the secondary integration.

\section{Numerical experiments}
\label{s:ne}

We next report some numerical experiments to asses our implementation of the $6$-stage Gauss collocation method of order $12$ in the 64-bit IEEE double precision  floating point arithmetic.

\subsection{Test problems}

We consider two different Hamiltonian problems corresponding to a double pendulum and the simulation of the outer solar system (considered in~\cite{Hairer2006} and \cite{Hairer2008}) respectively. In all the cases, we consider a time-step $h$ that is small enough for round-off errors to dominate over truncation errors. 

\subsubsection{The double pendulum (DP) problem}
 We consider the planar double pendulum problem: a double bob pendulum with masses $m_1$ and $m_2$ attached by rigid massless rods of lengths $l_1$ and $l_2$. 
This is a non-separable Hamiltonian system with two degrees of freedom, for which no explicit symplectic Runge-Kutta-type method is available, and hence Gauss collocation methods are a natural choice~\cite{McLachlan1992}.

\begin{figure} [h!]
\centerline{\includegraphics [width=10cm, height=8cm] {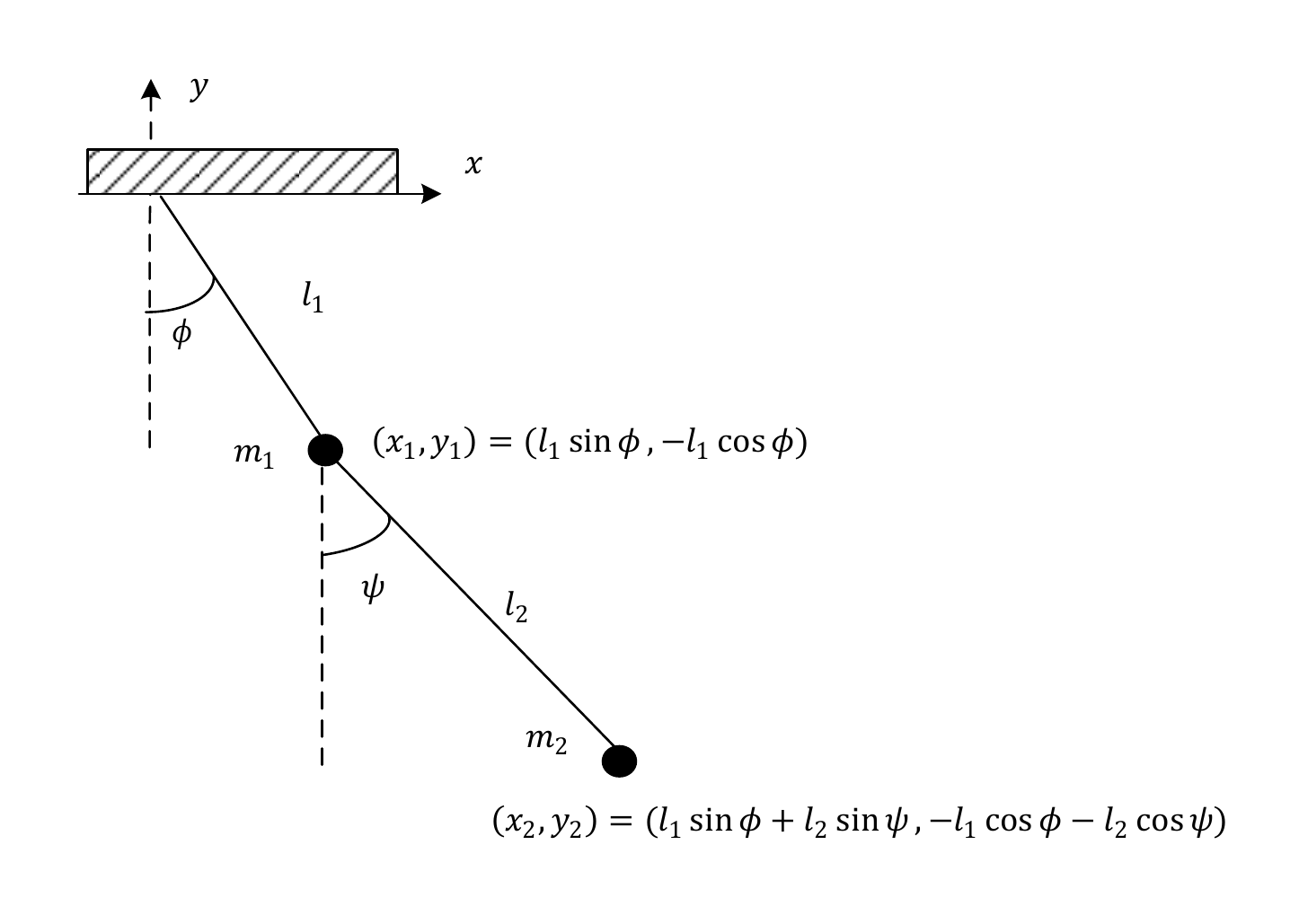}}
\caption{\label{fig:double-pendulum}Double Pendulum.}
\label{fig:two}
\end{figure}

{The configuration of the pendulum is described by two angles $q=(\phi,\theta)$  (see figure \ref{fig:double-pendulum}): while $\phi$ is the angle of the first bob, the second bob's angle is defined by $\psi=\phi+\theta$. We denote the corresponding momenta as $p=(p_{\phi},p_{\theta})$}.

Its Hamiltonian function $H(q,p)$ is
\begin{multline} \label{eq:2}
-\frac{ {l_1}^2 \ (m_1+m_2) \ {p_{\theta}}^2 +{l_2}^2 \ m_2 \ (p_{\theta} -p_{\phi})^2 + 2 \ l_1 \ l_2 \ m_2 \ p_{\theta} \ (p_{\theta} -p_{\phi}) \  \cos(\theta )} {{l_1}^2  \ {l_2}^2 \ m_2 \  (-2 \ m_1 - m_2 + m_2 \ \cos(2 \theta ))} \\
-g  \ \cos (\phi) \  (l_1 \ (m_1+m_2)+l_2 \ m_2 \ \cos(\theta))+g \ l_2 \ m_2 \ \sin(\theta) \sin(\phi),
\end{multline}
and we consider following parameters values
\begin{equation*} \label{eq:17}
g=9.8 \ \frac{m}{sec^2}\ ,\ \ l_1=1.0 \ m \ , \ l_2=1.0 \ m\ , \ m_1=1.0 \ kg\ , \ m_2=1.0 \ kg.
\end{equation*} 
 We take two initial values from \cite{Dumitru}, the first one of non-chaotic character, and the second one exhibiting chaotic behaviour:
\begin{enumerate}
   \item {Non-Chaotic case (NCDP): $q(0)=(1.1, -1.1)$ and $p(0)=(2.7746,2.7746)$.   We have integrated over $T_{end}=2^{12}$ seconds with step-size $h = 2^{-7}$. The numerical results will be sampled once every $m=2^{10}$ steps}.
   \item {Chaotic case (CDP): $q(0)=(0,0)$ and $p(0)=(3.873,3.873)$.
   We have integrated over $T_{end}=2^{8}$ seconds with step-size $h = 2^{-7}$.   We sample the numerical results  once every $m=2^{8}$ steps}.
\end{enumerate}
 Both initial value problems (NCDP and CDP) will be used to test the evolution of the global errors as well as to check the performance of the round-off error estimators. For the long term evolution of the errors in energy, only the NCDP problem will be considered.

\subsubsection{Simulation of the outer solar system (OSS)}

We consider a simplified model of the outer solar system (sun, the four outer planets, and Pluto) under mutual gravitational (non-relativistic) interactions.  This is a Hamiltonian system with 18-degrees of freedom ($q_i, p_i \in \mathbb{R}^3, \ i=0,\cdots,5$) and Hamiltonian function is
\begin{equation}
\label{eq:Ham2}
H(q,p)=\frac{1}{2}\ \sum^N_{i=0}{\ \frac{{\|p_i\|}^2}{m_i}}- \ G \sum^N_{0\le i<j\le N}{\frac{m_im_j}{\|q_i-q_j\|}}.
\end{equation}
We have considered the initial values and the values of the constant parameters ($G m_i$, $i=0,\ldots,5$) taken from ~\cite[page~14]{Hairer2006}. We  have integrated over $T_{end}=10^{7}$~days, with step-size $h=500/3$ and  the numerical results are sampled once every $m=120$ steps.

Observe that (\ref{eq:Ham2}) is a separable Hamiltonian, i.e.,  of the form $H(p,q) = T(p) + U(q)$.
It is well known that the efficiency of standard fixed point iteration can be improved for Hamiltonians systems with separable Hamiltonian  by considering a partitioned version of the fixed point iteration~\cite{Hairer2006}. Nevertheless, as in~\cite{Hairer2008}, here we report the numerical results obtained with the standard non-partitioned fixed-point iteration. (We have actually checked that similar results are obtained with the partitioned version of the fixed point iteration, the main difference being that less iterations are performed at each step.)

\subsection{Comparison of different sources of error in energy}

The error of a numerical solution $\tilde y_{n} + e_{n} \approx y(t_n)$  ($n=1,2,\ldots$) computed with our  double precision (DP) implementation of symplectic IRK schemes is a combined result of different kinds of errors:

\begin{enumerate}
\item The truncation error: The error due to replacing $y(t_n)$, $n=1,2,3,\ldots$ (where $y(t)$ is the solution of the initial value problem (\ref{eq:ivp})) by the numerical approximations  $y_{n}$ defined by (\ref{eq:y})--(\ref{eq:YL})  (with exact coefficients $b_{i}, \mu_{i j}$). 

\item The iteration error:  In practice a finite number $K$ of fixed point iterations (\ref{eq:fixed_point_iteration}) are applied, and the solution $L_{n,i}, Y_{n,i}$ ($i=1,\ldots,s$) of (\ref{eq:YL}) are replaced by approximations  $L_{n,i}^{[K]}, Y_{n,i}^{[K]}$. 
Then, in the FPIEA implementation, the corresponding numerical solution $\overline{y}_{n+1}$ is computed at each step as
\begin{equation*}
\overline{y}_{n+1} = y_{n} + \sum_{i=1}^{s} L_{n,i}^{[K]}.
\end{equation*}

\item The error due to replacing the original map $f:\R^D \to \R^D$ by its computational substitute 
$\tilde{f}:\F^D \to \F^D$.  This has a double effect: From one hand, in most steps, a computational fixed point is achieved in a finite number $K$ of iterations, which {causes} an unavoidable iteration error. On the other hand, replacing $f$ by $\tilde f$ adds the effect of some round-off errors. 

\item The error due to the application of a different IRK scheme. In our case, we apply the scheme (\ref{eq:y})--(\ref{eq:YL}) with $b_i$ replaced by $\tilde b_i \in \F$ ($i=1,\ldots,s$) and each $\mu_{i j}$  replace by double precision approximation  $\tilde \mu_{i j} \in \F$ satisfying condition (\ref{eq:sympl_cond_2}).

\item The error due to using inexact arithmetic for the operations (other than the evaluation of $\tilde f$) in the machine precision implementation of the algorithm.
\end{enumerate}

We have simulated, for the double pendulum (the non-chaotic case,  NCDP) and the outer solar system (OSS) respectively, the effect that each of the first four of such sources of errors has in the values of the energy (which of course is conserved in the exact solution) as follows:

\begin{figure}[h!]
\centering
\begin{tabular}{c c}
\subfloat[NCDP ($h=2^{-7}$)]
{\includegraphics[width=.600\textwidth]{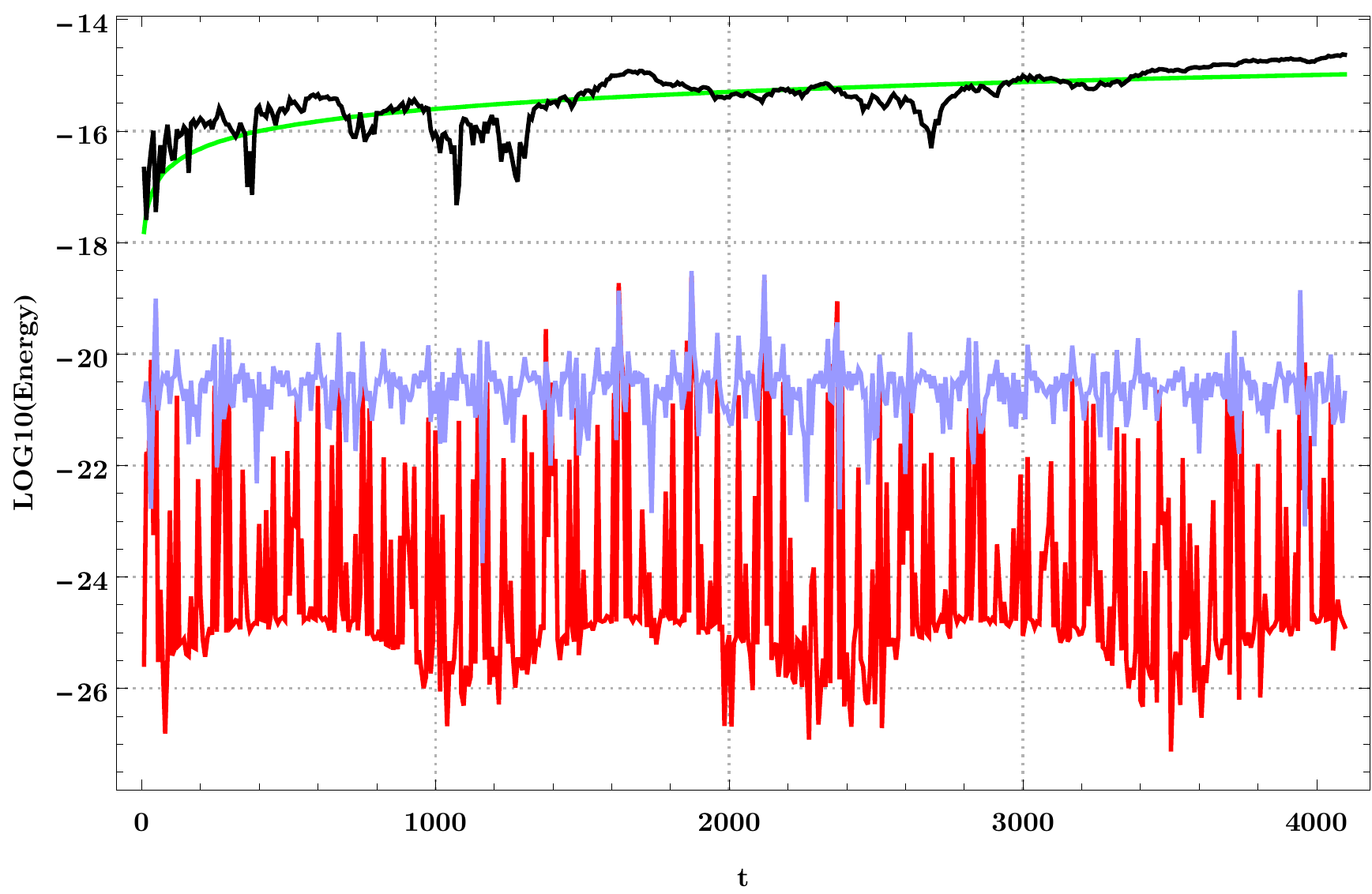}}\\
\subfloat[OSS ($h=500/3$)]
{\includegraphics[width=.600\textwidth]{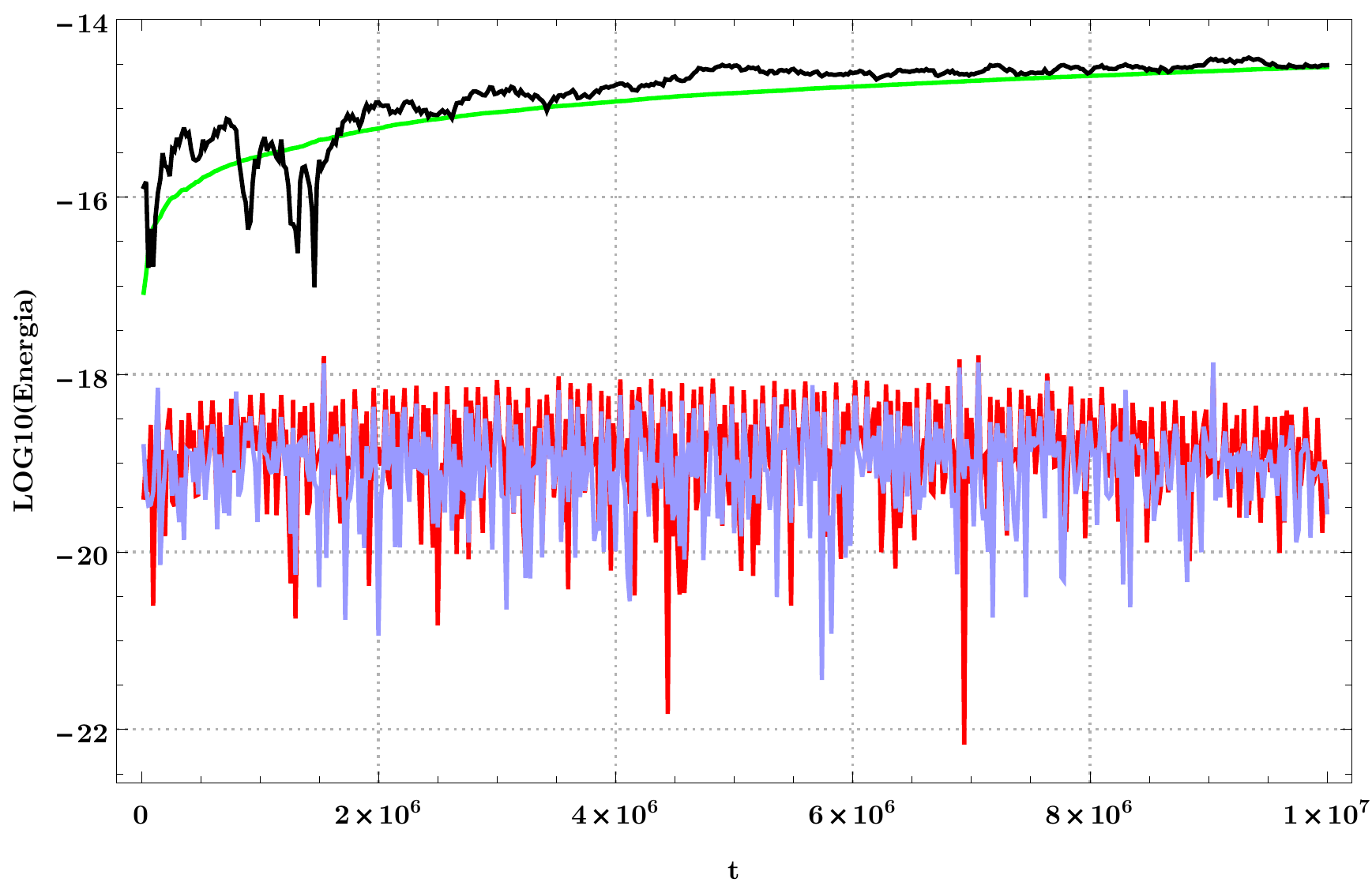}}
\end{tabular}
\caption{\small We plot the evolution of energy error in logarithmic scale of the next algorithms implementations: A-algorithm  as estimation of truncation error (red), B-algorithm  as estimation of iteration error (green), C-algorithm  as estimation of the effect of replacing the exact $f$ by its double precision version $\tilde{f}$ function (black) and D-algorithm  as estimation of the effect of using double precision coefficients (blue).}    
\label{fig:plot1}
\end{figure}

\begin{enumerate}
\renewcommand{\theenumi}{\Alph{enumi}}
\item In order to estimate the truncation error,  we have applied our algorithm fully implemented in quadruple precision.
\item For the iteration error, we have applied the quadruple precision version of the algorithm modified so that the fixed point process in the $n$th step is stopped at the $K$th iteration provided that $Y_{n,i}^{[K]}$ and $Y_{n,i}^{[K-1]}$ coincide when rounded to double precision. 
\item In addition, we have estimated the  effect (in the evolution of the energy) of the error due to replacing $f$ by $\tilde f$, by considering the quadruple precision implementation of our algorithm with the double precision version of $\tilde f$.
\item Finally, we have simulated the error due to the application of a RK scheme with double precision coefficients 
 by applying our quadruple precision implementation with  double precision coefficients $\tilde b_i, \tilde \mu_{i j}$.
 \end{enumerate}

We next plot (Fig. \ref{fig:plot1}), for each of the considered initial value problems, the evolutions of the energy errors corresponding to the items A--D in previous list. In both cases, we have chosen a step-size $h$ such that truncation errors are smaller than round-off errors.  We observe that the effect of using  double precision coefficients ($\tilde{b}_i, \tilde{\mu}_{ij}$) is also negligible compared to the propagation of round-off errors. The iteration error is  similar in size to round-off errors.

\subsection{Statistical analysis of errors}

In order to make a more robust comparison of the numerical errors due to round-off errors, we adopt  (as in~\cite{Hairer2008}) an statistical approach. We have considered for each of the three initial value problem, $P=1000$ perturbed initial values by randomly perturbing each component of the initial values with a relative error of size $\mathcal{O}(10^{-6})$.

We will compare three different fixed point implementations of the $6$-stage Gauss collocation method. In all of them, the same computational substitute $\tilde f:\F^D \to \F^D$ is used instead of the original map $f:\R^D \to \R^D$ defining the ODE (\ref{eq:ivp}):

\begin{enumerate}
\item  The FPIEA (fixed point iteration with exact arithmetic) implementation, where the techniques described in Subsections ~\ref{ss:3.1} and~\ref{ss:3.2} are applied to implement a fixed point iteration with all arithmetic operations (other than those used when evaluating $\tilde f$) performed in exact arithmetic. 
\item Our double precision version (coded in C) of the algorithm implemented in FPIEA.
 We will refer to it as DP (double precision).
\item The algorithm proposed in~\cite{Hairer2008}, implemented in \href{http://www.unige.ch/~hairer/preprints/code.tar}{Hairer's Fortran code}.
\end{enumerate}
From one hand, we want to check if the propagation of round-off errors in our DP implementation are  qualitatively similar and close in magnitude to its exact arithmetic counterpart FPIEA. On the other hand, we want to see how our DP implementation compares with Hairer's code.

In (Table~\ref{tab:fp}) we display  the percentage of steps that reach a computational fixed-point and the average number of fixed-point iterations per step in each of the referred three implementations when applied to two differenent initial value problems.

\begin{table}
\caption[Fixed-point percentage of steps and mean iterations.] 
{\small{Percentage of steps that reach a computational fixed-point and the number of fixed-point iterations per step for the computations of non-chaotic double pendulum (NCDP), chaotic double pendulum (CDP), and the outer solar system (OSS) problems. In columns, we compare three different implementations: FPIEA, DP (double precision) and Hairer's Fortran code.}}
\label{tab:fp}       % Give a unique label
\centering
%\resizebox{\textwidth}{!}
{%
\begin{tabular}{ l c c c c c c } 
 \hline
                 &  \multicolumn{2}{c}{FPIEA}  & \multicolumn{2}{c}{DP} & \multicolumn{2}{c}{Hairer} \\
 \hline
 NCDP            & $98.7\%$    & $8.5$   & $98.8\%$     & $8.6$   &  $98.5\%$ & $8.6$  \\ 
 CDP             & $98.9\%$    & $8.5$   & $98.9\%$     & $8.6$   &  $98.4\%$ & $8.6$  \\ 
 OSS             & $97.7\%$    & $14.4$  & $97.4\%$     & $14.2$  &  $87.5\%$ & $14.1$ \\ 
   \hline
 \end{tabular}}
\end{table}

\subsubsection{Distribution of energy jumps}

\begin{figure}[h!]
\centering
\begin{tabular}{c c}
\subfloat[\small {NCDP}]
{\includegraphics[width=.4\textwidth]{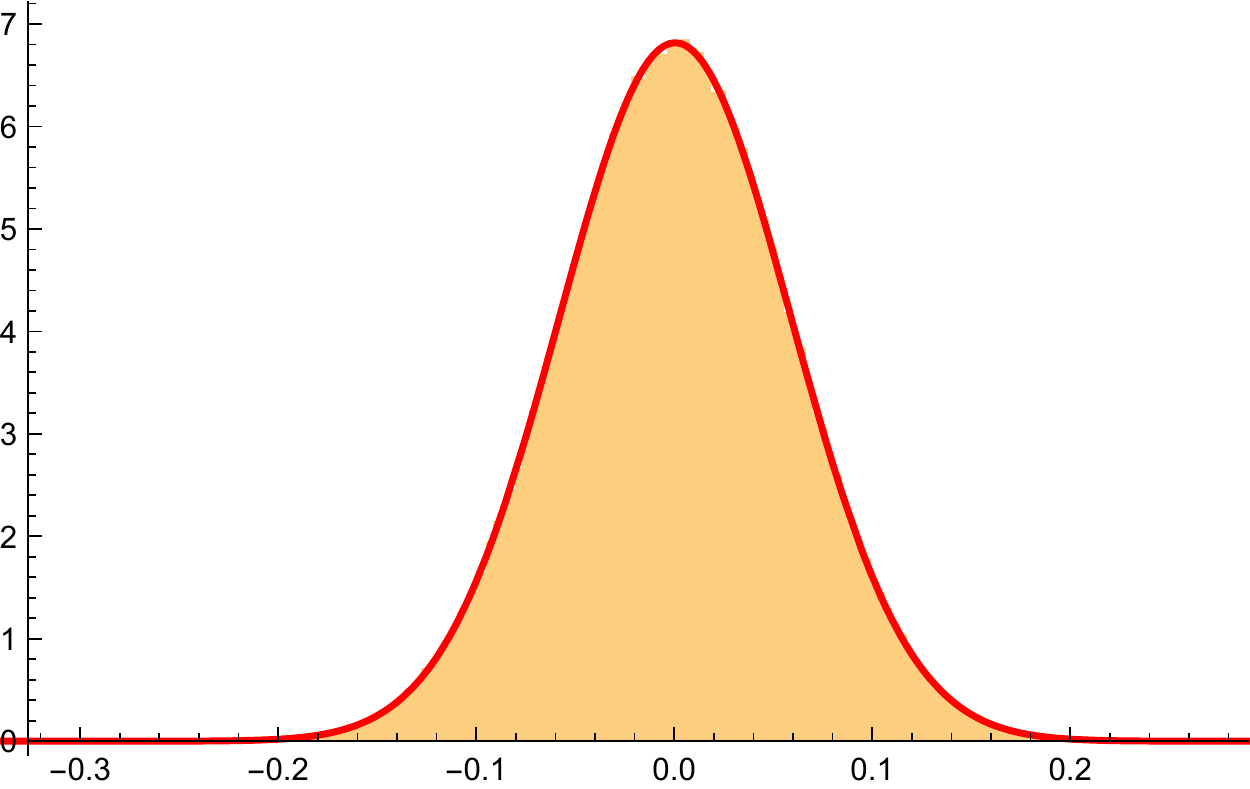}}
&
\subfloat[OSS]
{\includegraphics[width=.4\textwidth]{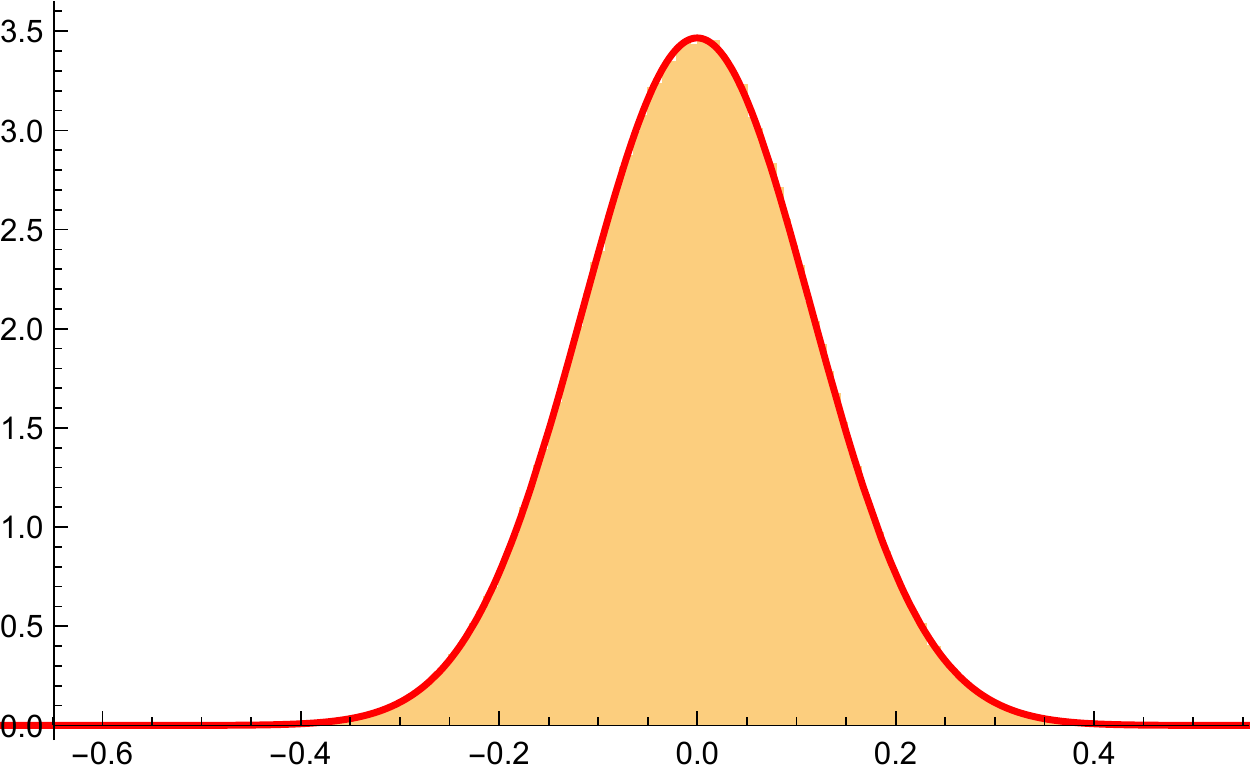}}
\\
 ($\mu=5.3\times 10^{-19}$, \ $\sigma=1.5\times 10^{-17}$) &
 ($\mu=-1.9\times 10^{-19}$, \ $\sigma=3.5\times 10^{-18}$) 
 \end{tabular}
\caption{ \small Histograms of $K P$ samples of energy jumps of the  DP implementation against the normal distribution $N(\mu, \delta)$ for two problems (NCDP and OSS). The horizontal axis is multiplied by $10^{15}$ and vertical axis indicates the frequency }
\label{fig:plot2}
\end{figure}

The local error in energy $H(y_n)-H(y_{n-1})$ due to round-off errors, is  "expected" to behave, for a good implementation free from biased errors,  like an independent random variable.  Then, provided that the numerical results are sampled every $m$ steps, 
with a large enough sampling frequency $m$, an energy jump $H(y_{k m})-H(y_{mk-m})$ will behave as an independent variable with an approximately Gaussian distribution with mean $\mu$ (ideally $\mu=0$) and standard deviation $\sigma$.  So that the accumulated difference in energy, 
\begin{equation}
\label{eq:AEE}
H(y_{k m})-H(y_{0})
\end{equation}
at the sampled times $t_{m k} = k m h$ would  behave like a Gaussian random walk with standard deviation  $k^{\frac12}\sigma = (t_{m k}/(m h))^{1/2} \sigma$.  This is sometimes referred to as Brouwer's law in the scientific literature~\cite{Grazier2005}, from the original work on the accumulation of round-off errors in the numerical integration of Kepler's problem done by Brouwer in~\cite{Brouwer1937}.

In this sense, we want to check in which extent the (scaled) energy jumps,
\begin{equation}
\label{eq:REJ}
(H(y_{k m})-H(y_{mk-m}))/H(y_0)
\end{equation}
 due to round-off errors after $m$ steps approximately obey a Gaussian distribution in our double precision (DP) implementation. 
 
 If $[0,T_{end}]$ is the integration interval, and $P$ perturbed initial values are considered, we have a total number of $K P$ samples of energy jumps, where $K=T_{end}/(m h)$. In Figure~\ref{fig:plot2}, we plot the histograms of $K P$ samples of energy jumps obtained with our DP implementation against the normal distribution $N(\mu, \delta)$ (where $\mu$ and $\sigma$ are the average and standard deviation of the samples respectively). For both initial value problems, non-chaotic double pendulum (NCDP), and the outer solar system (OSS), such histograms fits perfectly well to their corresponding normal distributions $N(\mu, \delta)$.

\subsubsection{Evolution of mean and standard deviation of errors}

\begin{figure}[h!]
\centering
\begin{tabular}{c c}
\subfloat[NCDP: energy error]
{\includegraphics[width=.4\textwidth]{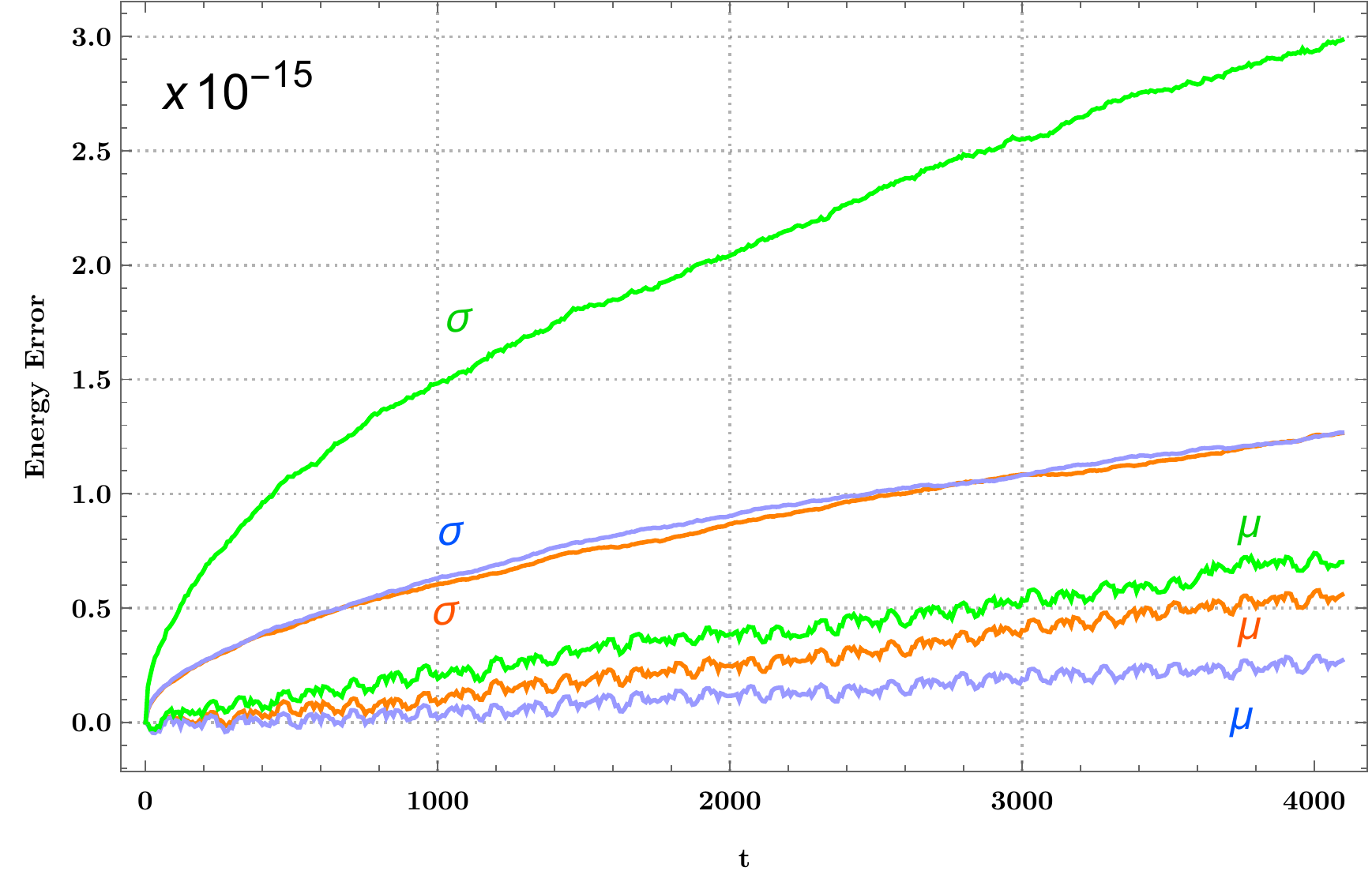}}
&
\subfloat[NCDP: detail of mean error.]
{\includegraphics[width=.4\textwidth]{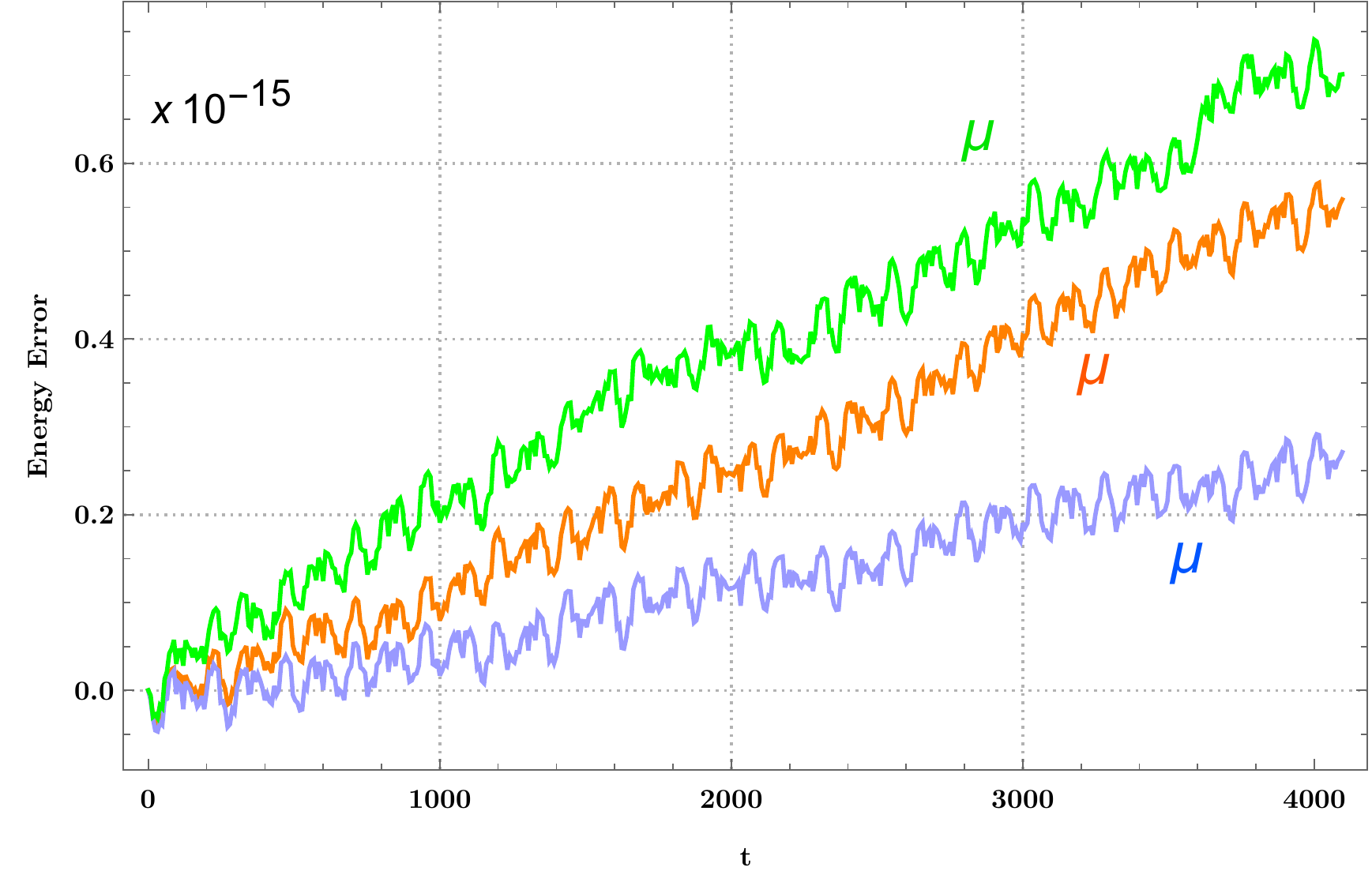}}
\\
\subfloat[OSS: energy error.]
{\includegraphics[width=.4\textwidth]{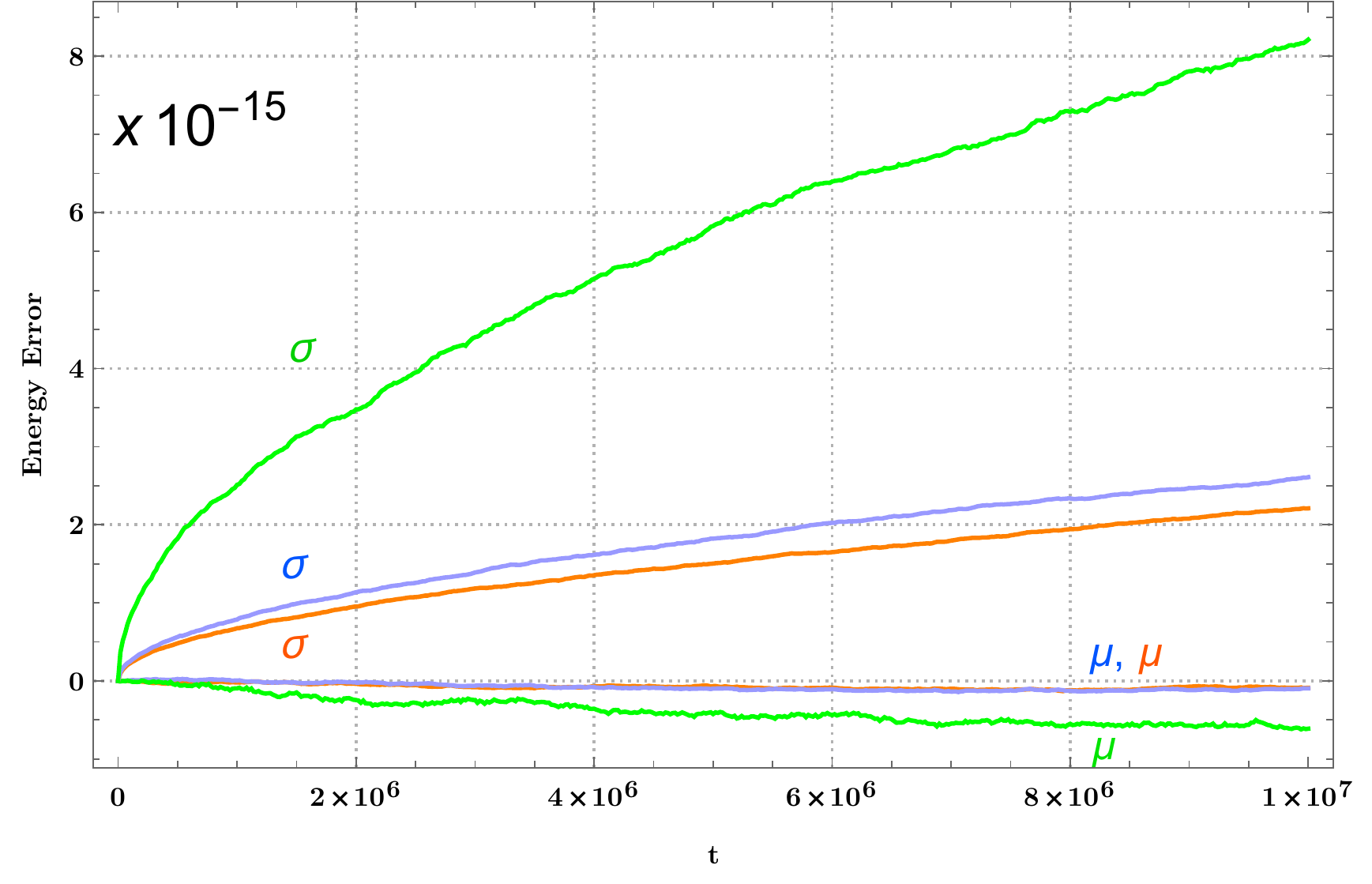}}
&
\subfloat[OSS: detail of mean error.]
{\includegraphics[width=.4\textwidth]{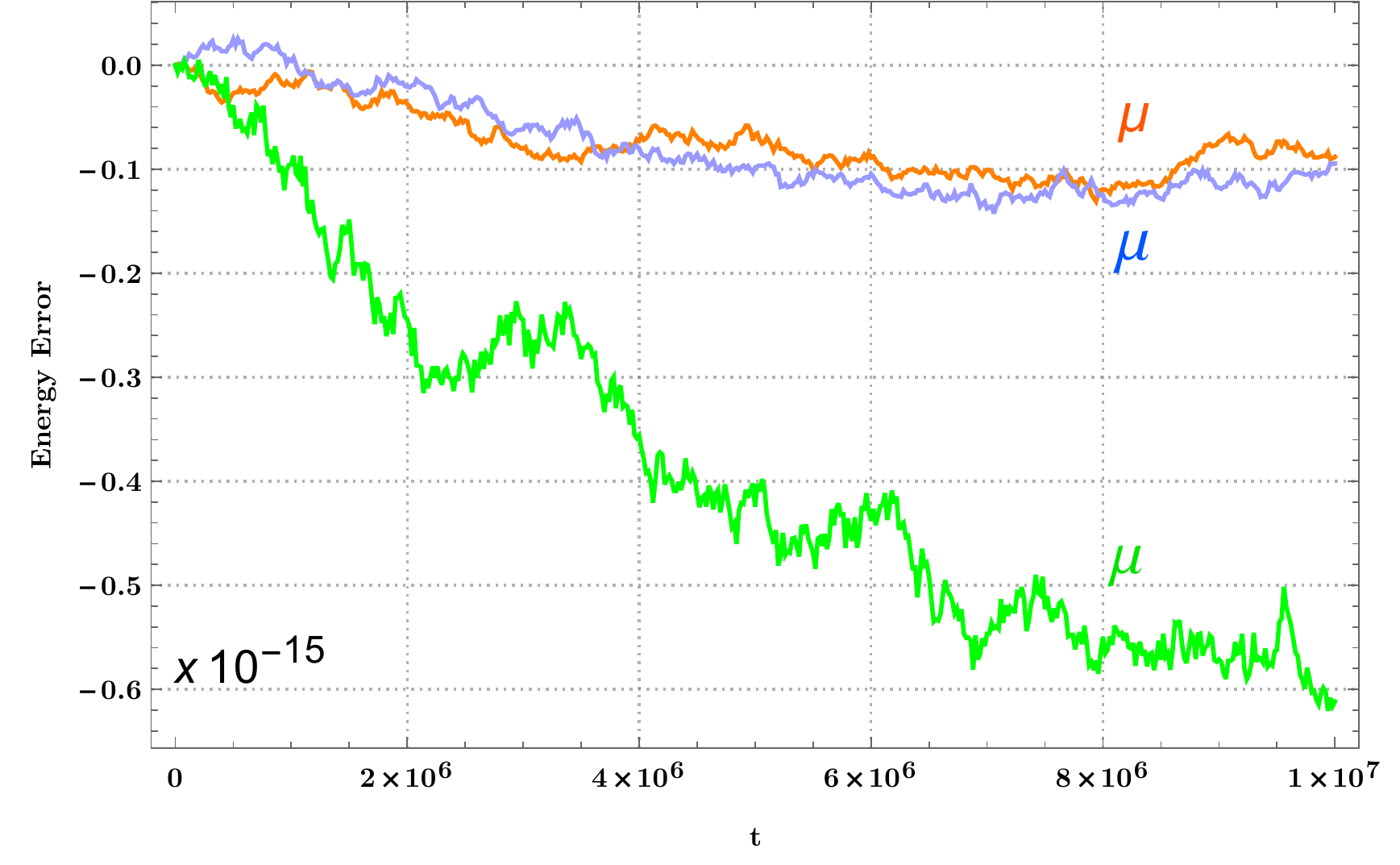}}
\end{tabular}
\caption{\small Evolution of mean ($\mu$) and standard deviation ($\sigma$) of errors in energy (left) and detail of the evolution of mean errors in energy (right), for DP implementation (blue),  FPIEA implementation (orange), and  Hairer's implementation (green).  Non-Chaotic case (a,b), and outer solar system case (c, d)}
\label{fig:plot3}
\end{figure}

We next plot (Fig.~\ref{fig:plot3}) the evolution of the mean and standard deviation of the errors in energy  of the FPIEA, DP, and Hairer's implementations  for the NCDP and OSS problems respectively.

Recall that FPIEA represents the best possible fixed point implementation of the IRK scheme provided that the double precision version $\tilde f$ of the original $f$ is used. We stress that we have made the stopping criterion of the FPIEA implementation even more stringent than in the DP implementation: we stop the fixed point iteration if  either $ \Delta^{[k]} =0$ or (\ref{eq:stopping}) is fulfilled during ten consecutive iterations.
This way, we try to avoid the persistence of iteration errors in the case of steps where no computational fixed point is obtained. 
(Observe that whenever $ \Delta^{[k]} =0$, there is no point in performing more fixed point iterations, as in that case a computational fixed point has been achieved.)

The numerical tests in Figure~\ref{fig:plot3} seen to confirm that our DP implementation is near optimal (that is, close to the FPIEA implementation), both with respect to the standard deviation and the mean of the errors in energy.

We believe that some small linear drift of the mean energy error may be unavoidable for the fixed point implementations of IRK schemes in some cases (such us the NCDP example). This is consistent with the observation that the simulated iteration error displaying in 
Figure~\ref{fig:plot0} is close in magnitude to the effect of round-off errors.  

\begin{figure}[h!]
\centering
\begin{tabular}{c c}
\subfloat[NCDP: energy error.]
{\includegraphics[width=.4\textwidth]{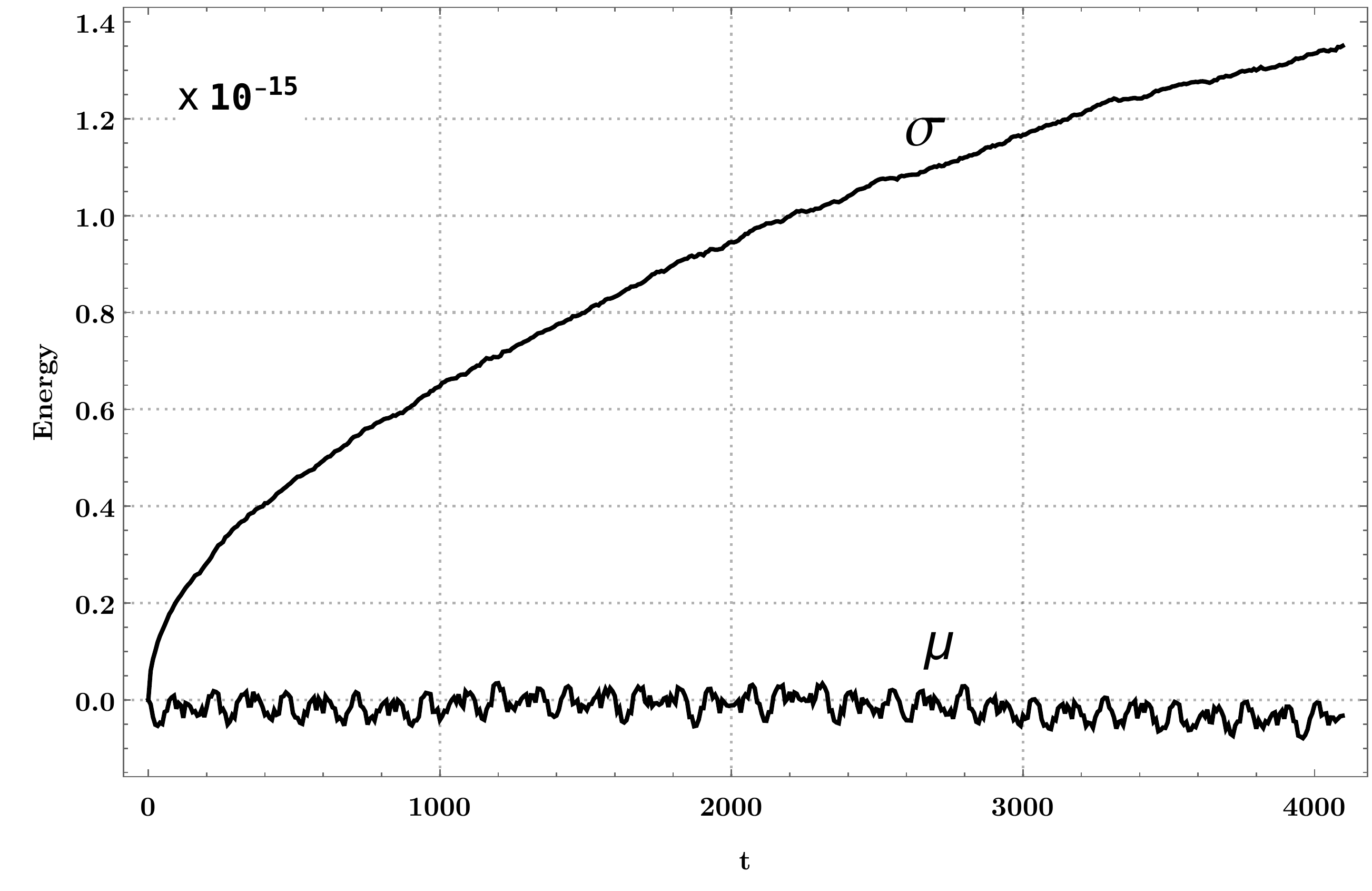}}
&
\subfloat[NCDP: detail of mean error.]
{\includegraphics[width=.4\textwidth]{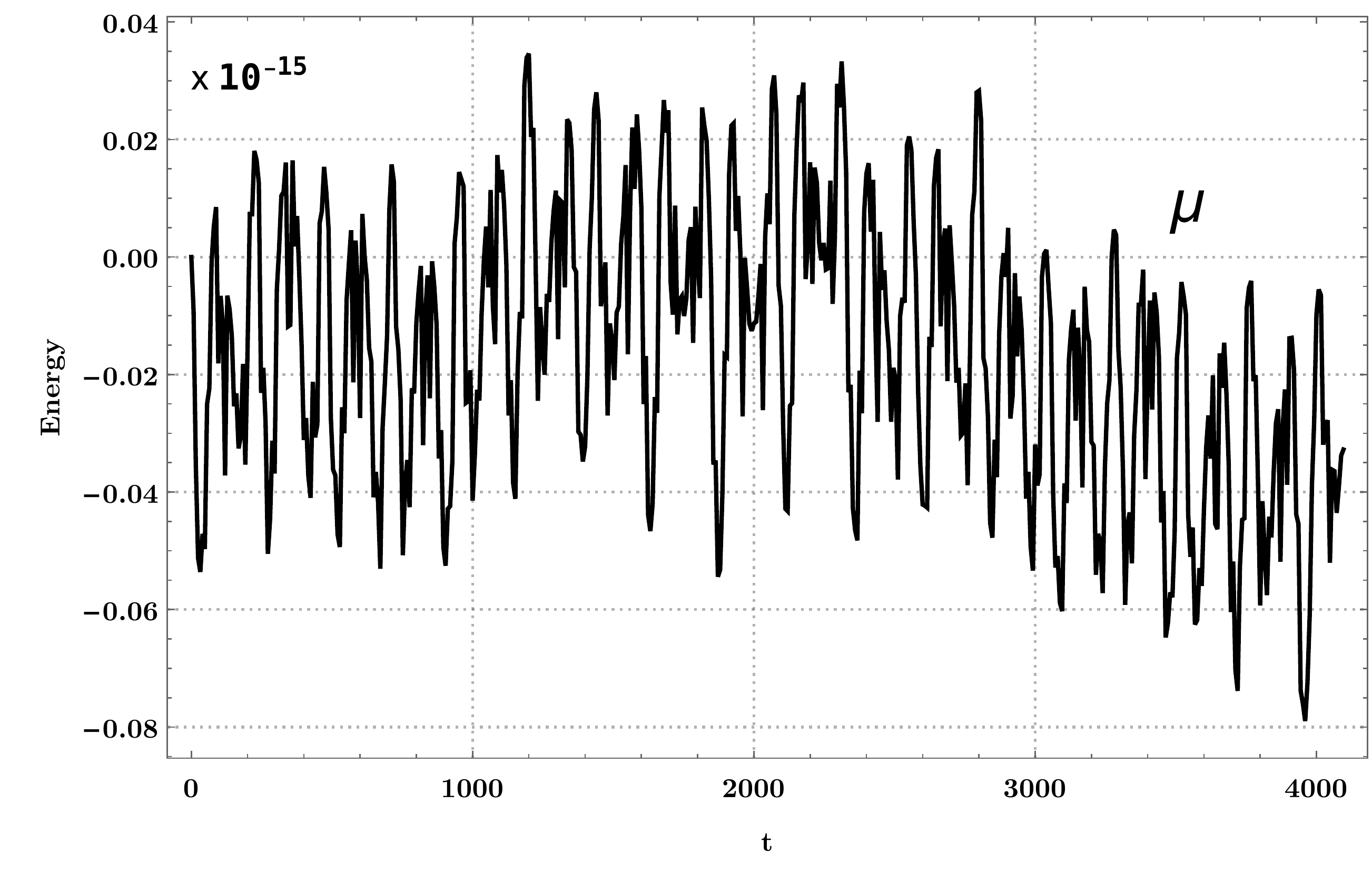}}
\end{tabular}
\caption{\small Evolution of mean ($\mu$) and standard deviation ($\sigma$) of errors in energy of a IRK implementation with simplified Newton iterations}
\label{fig:plotNewton}
\end{figure}

This is not of course inherent to the symplectic IRK scheme itself. In Figure~\ref{fig:plotNewton}, we display the results obtained for the NCDP example with a preliminary implementation of a simplified Newton implementation of the same IRK scheme as above.  No linear drift of the mean energy error seems to appear.

\begin{figure}[h!]
\centering
\begin{tabular}{c c}
\subfloat[NCDP: mean global error.]
{\includegraphics[width=.4\textwidth]{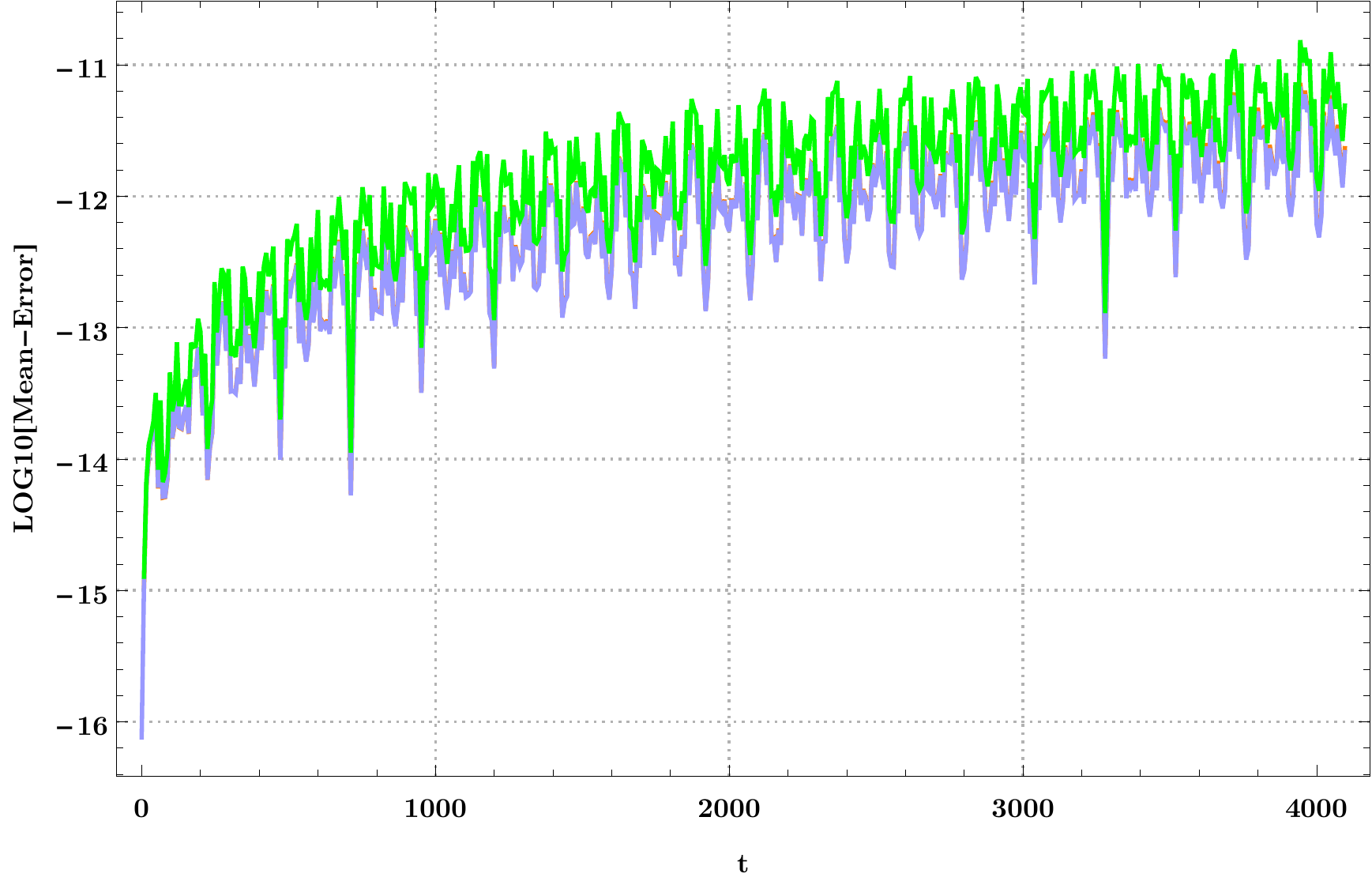}}
&
\subfloat[NCDP: standard deviation global error.]
{\includegraphics[width=.4\textwidth]{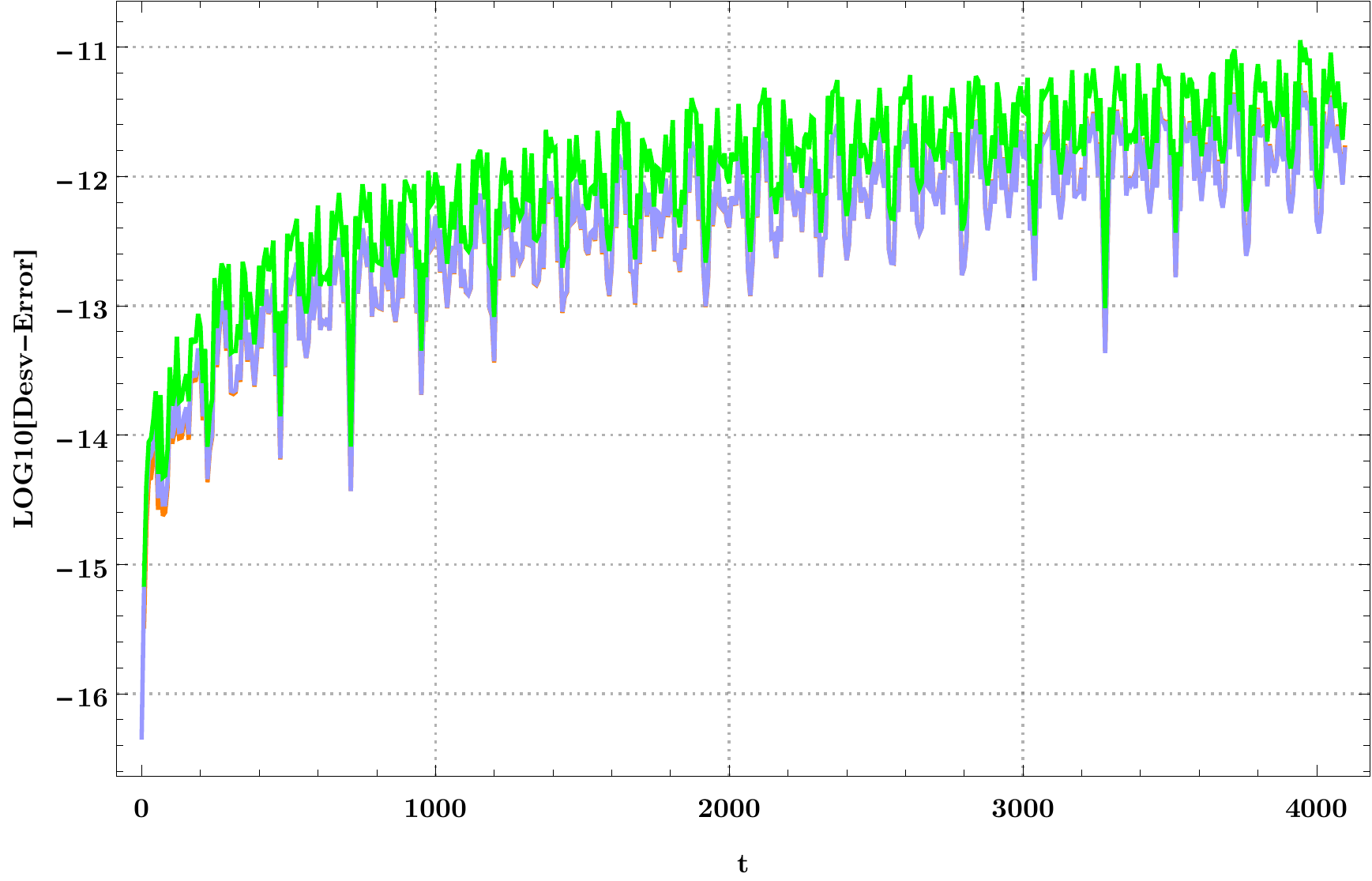}}
\\
\subfloat[CDP: mean global error.]
{\includegraphics[width=.4\textwidth]{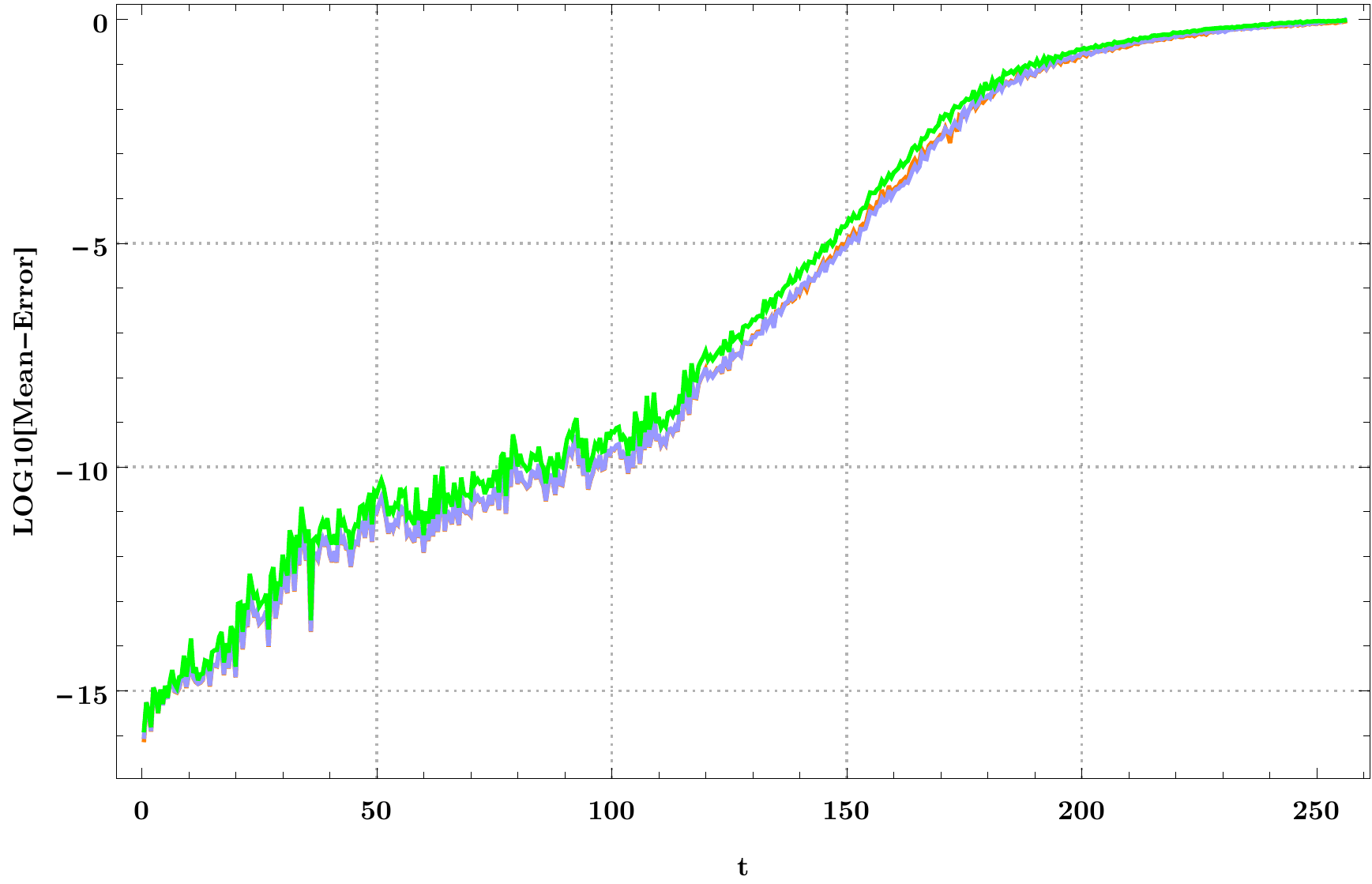}}
&
\subfloat[CDP: standard deviation global error.]
{\includegraphics[width=.4\textwidth]{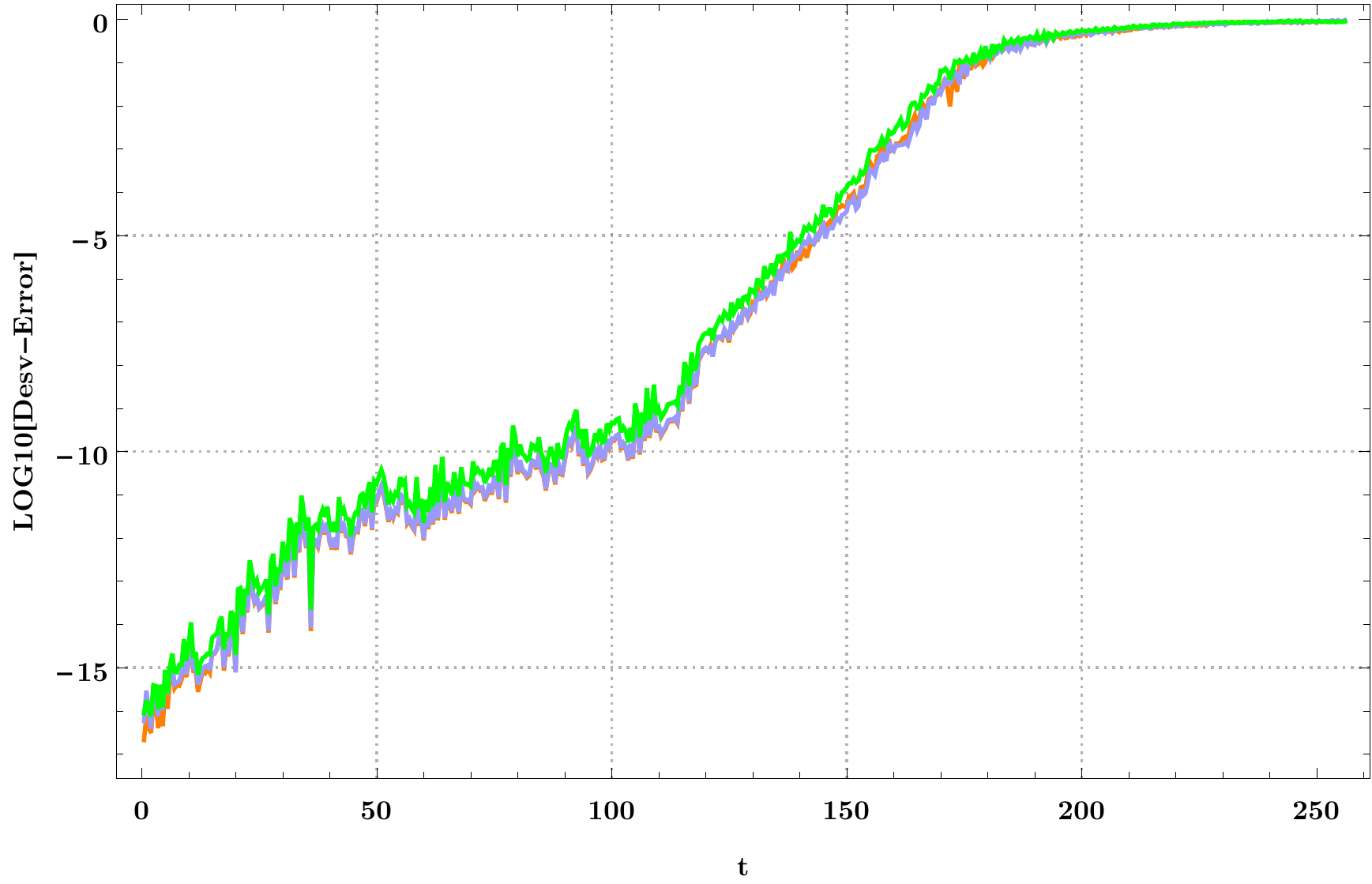}}
\\
\subfloat[OSS: mean global error.]
{\includegraphics[width=.4\textwidth]{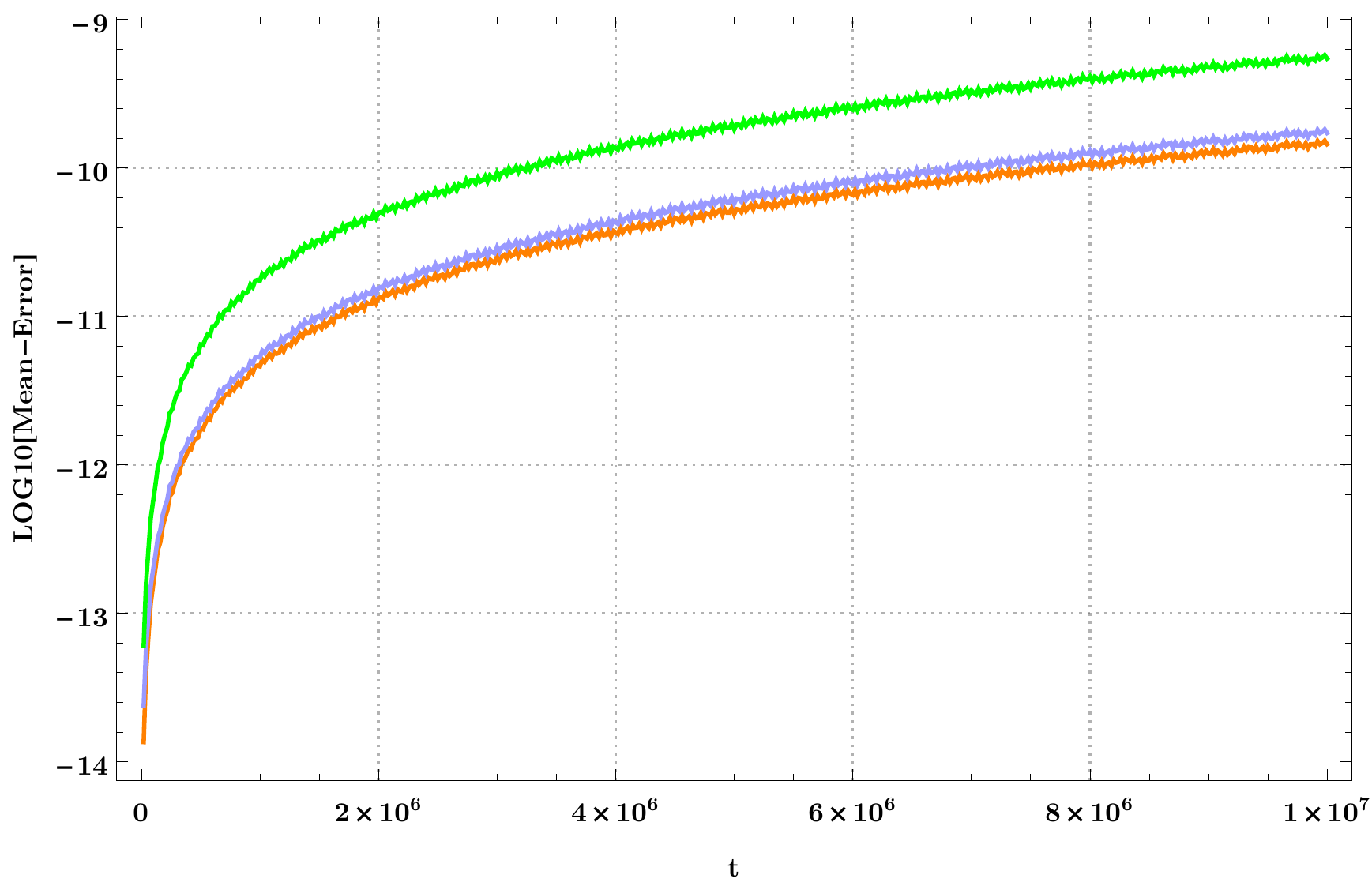}}
&
\subfloat[OSS: standard deviation global error.]
{\includegraphics[width=.4\textwidth]{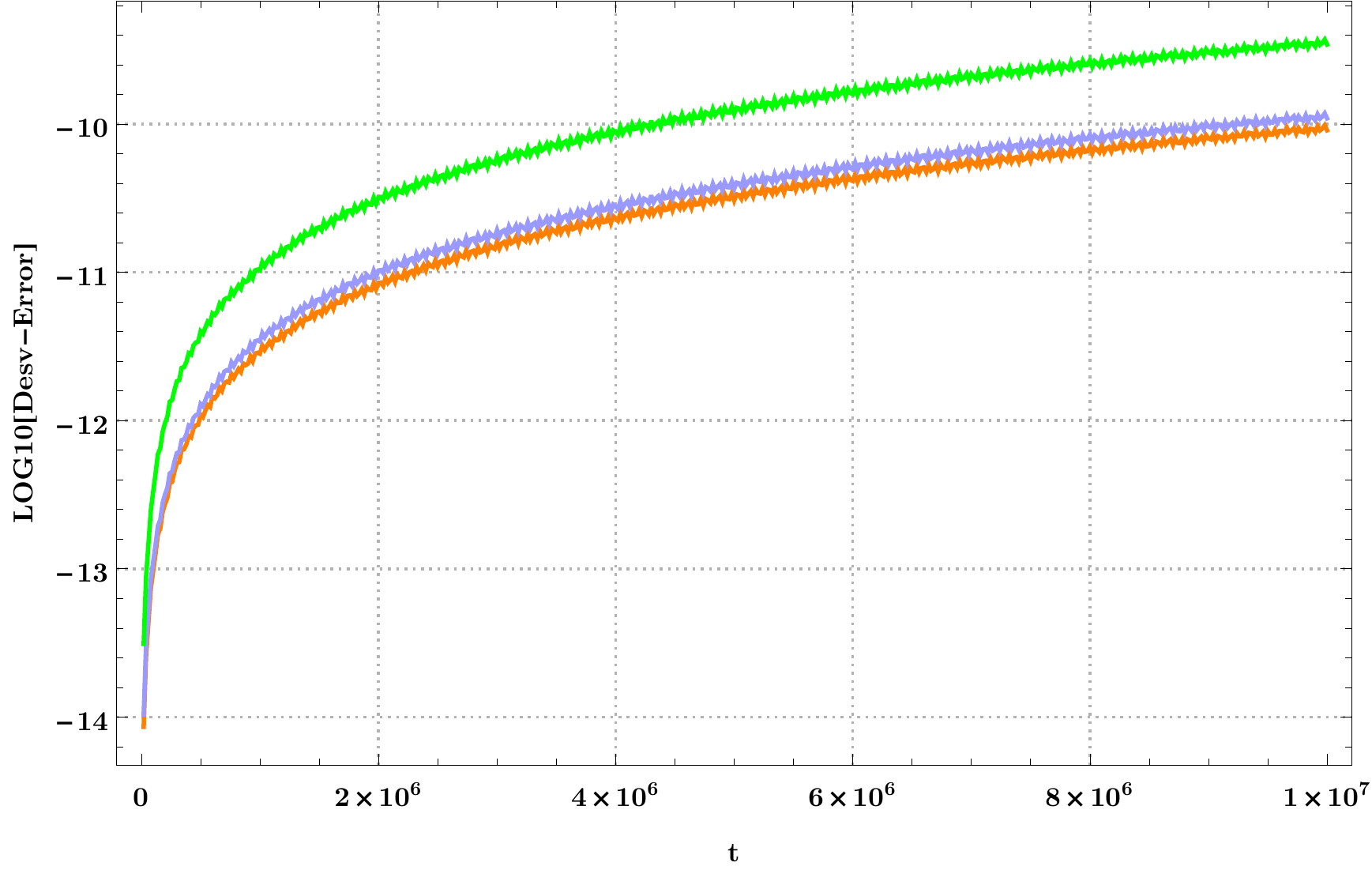}}
\end{tabular}
\caption{\small Evolutions of Mean (left)  and standard deviation (right) of global errors in positions of DP implementation (blue), FPIEA implementation (orange), and  Hairer's implementation (green):  NCDP (a,b), CDP (c,d) and OSS (e, f)}
\label{fig:plot4}
\end{figure}

To end with the present subsection, we plot (Fig.~\ref{fig:plot4}) the evolution of the (mean and standard deviation of) errors in position  of the FPIEA,  DP, and Hairer's implementations  for the NCDP, CDP and OSS problems respectively. The displayed results seem to confirm our claim of the DP implementation being a close-to-optimal fixed point implementation of the symplectic IRK scheme.

\subsection{Round-off error estimation}

\begin{figure}[h!]
\centering
\begin{tabular}{c c}
\subfloat[NCDP: original initial values]
{\includegraphics[width=.4\textwidth]{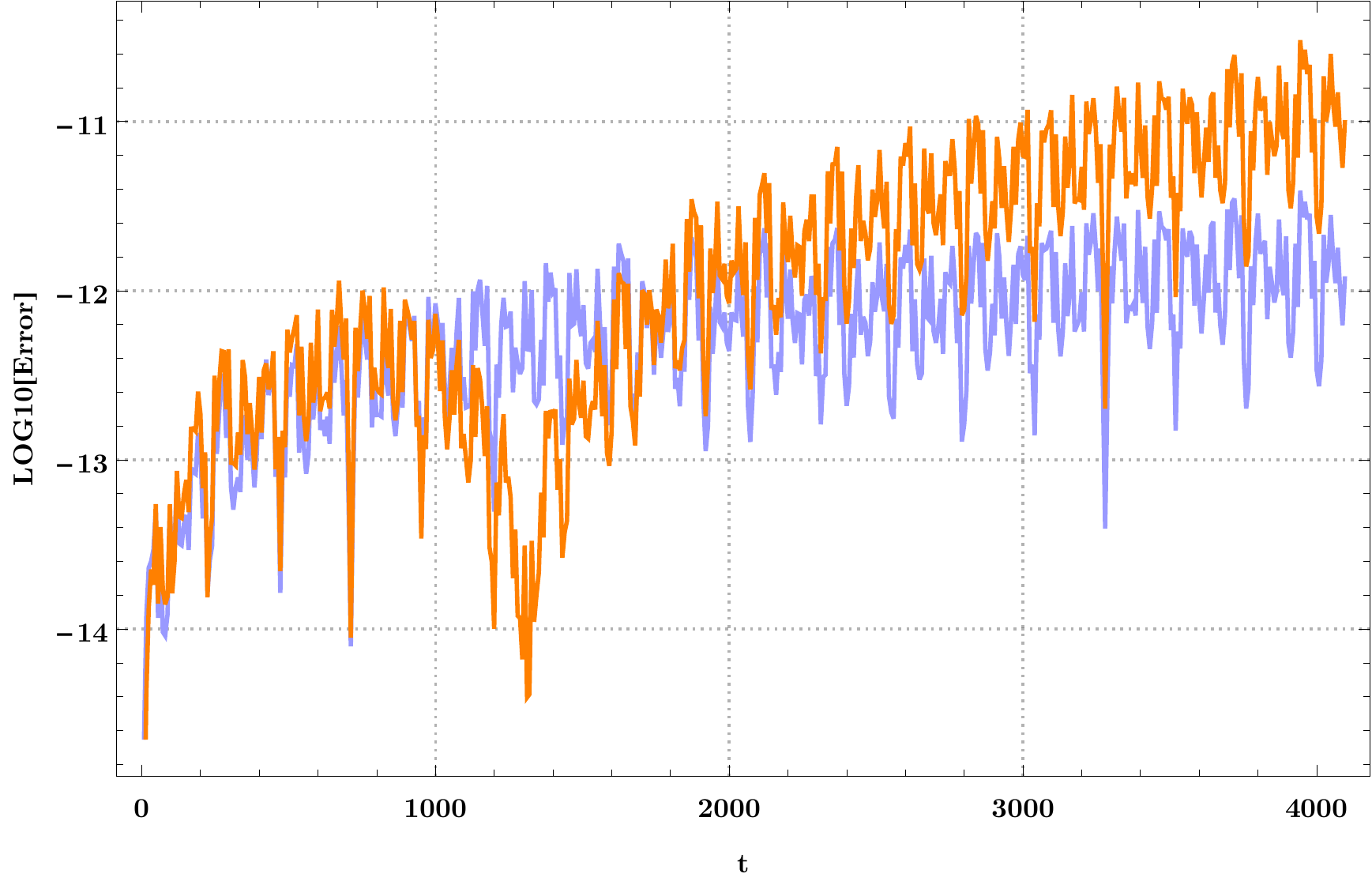}}
&
\subfloat[NCDP: $P=1000$ perturbed initial values]
{\includegraphics[width=.4\textwidth]{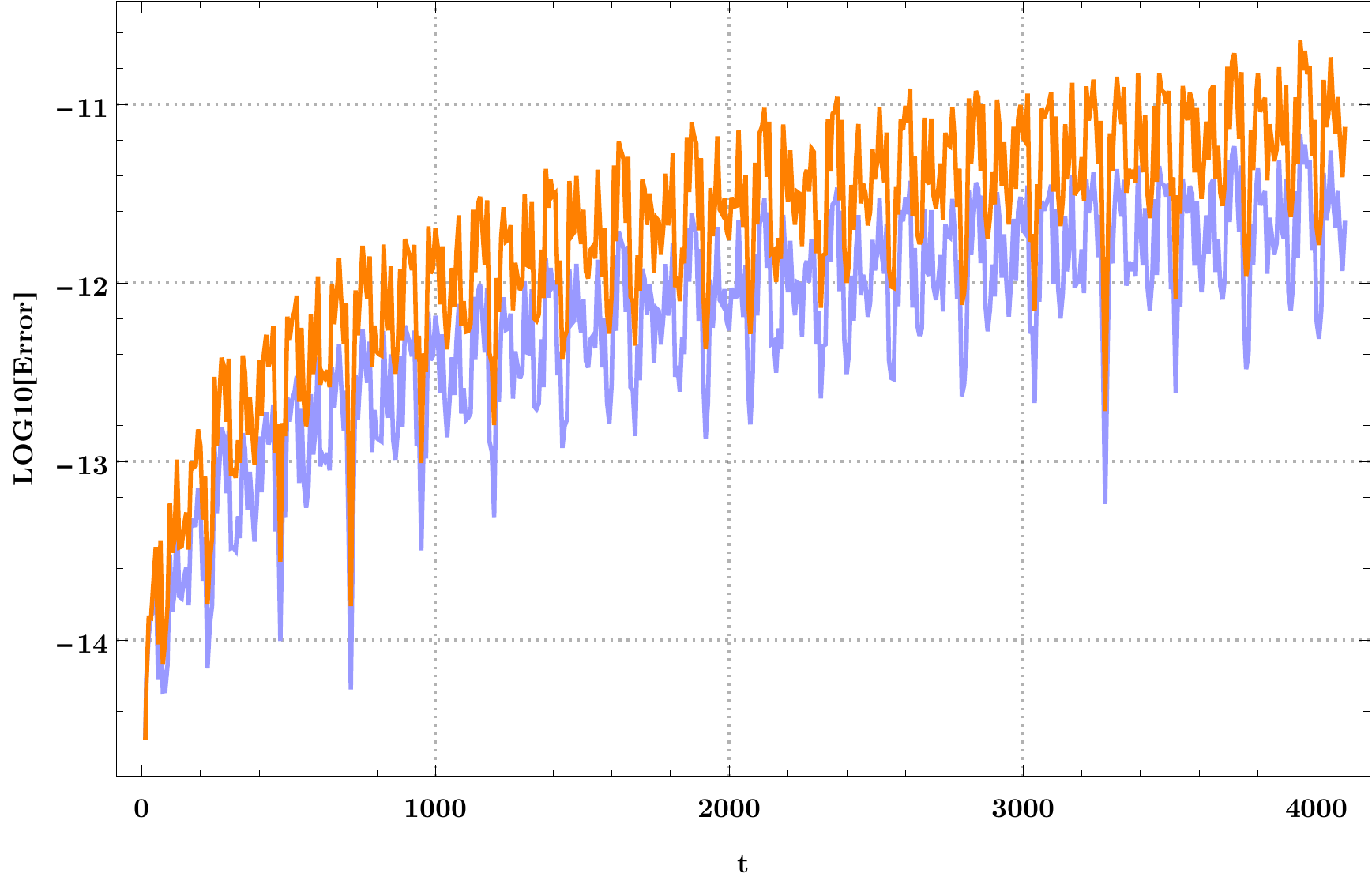}}
\\
\subfloat[CDP: original initial values]
{\includegraphics[width=.4\textwidth]{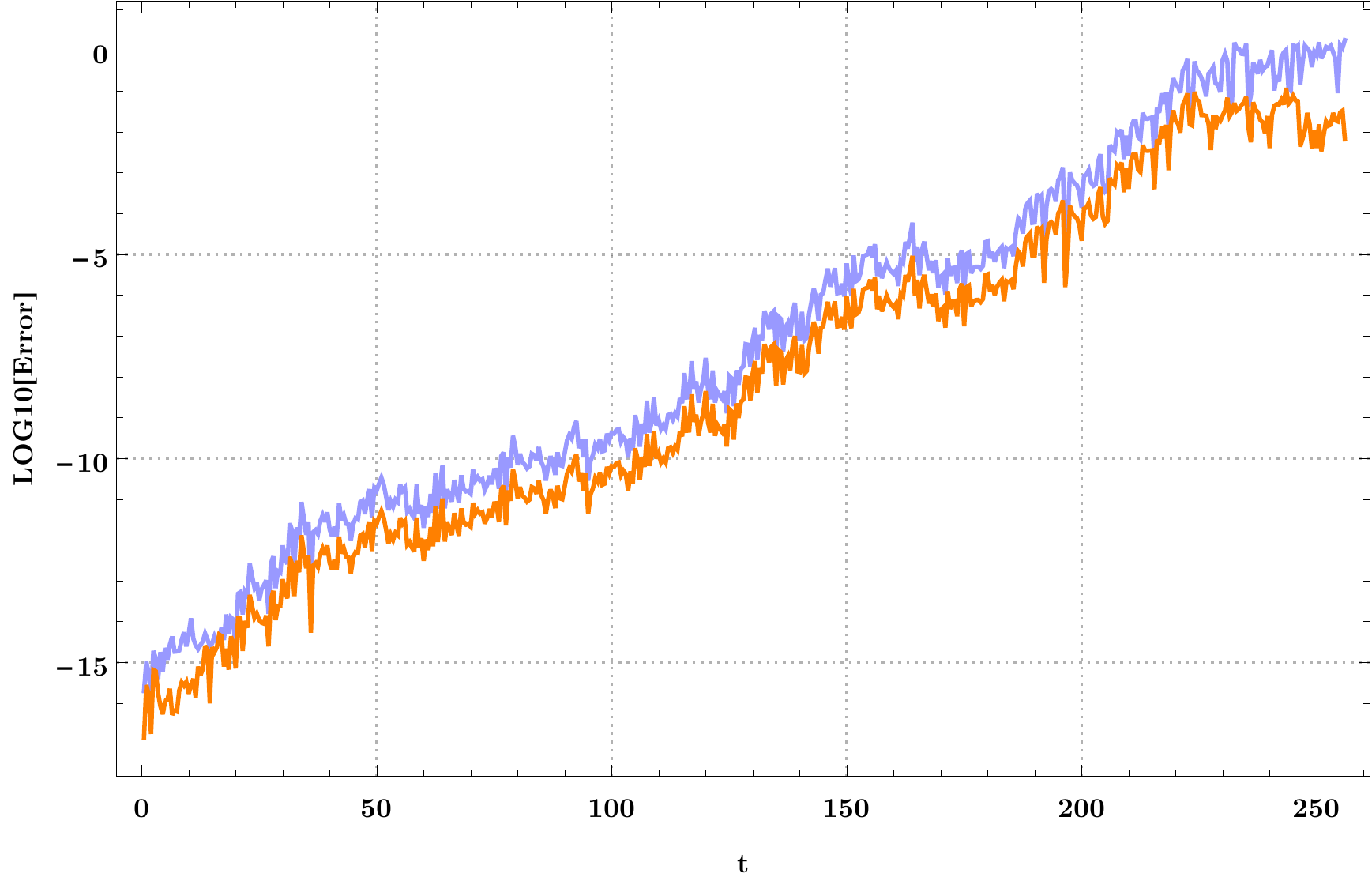}}
&
\subfloat[CDP: $P=1000$ perturbed initial values]
{\includegraphics[width=.4\textwidth]{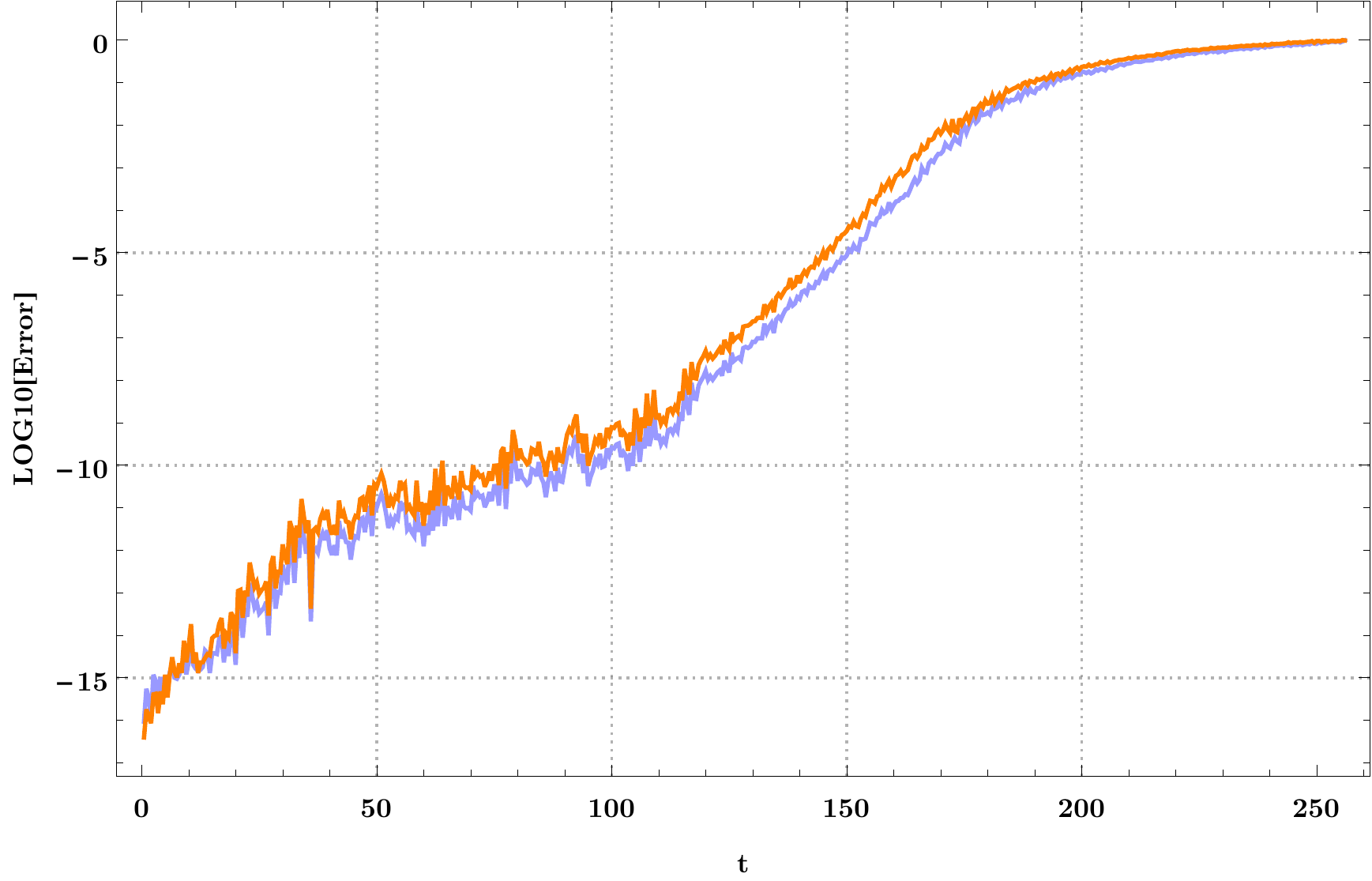}}
\\
\subfloat[OSS: original initial values]
{\includegraphics[width=.4\textwidth]{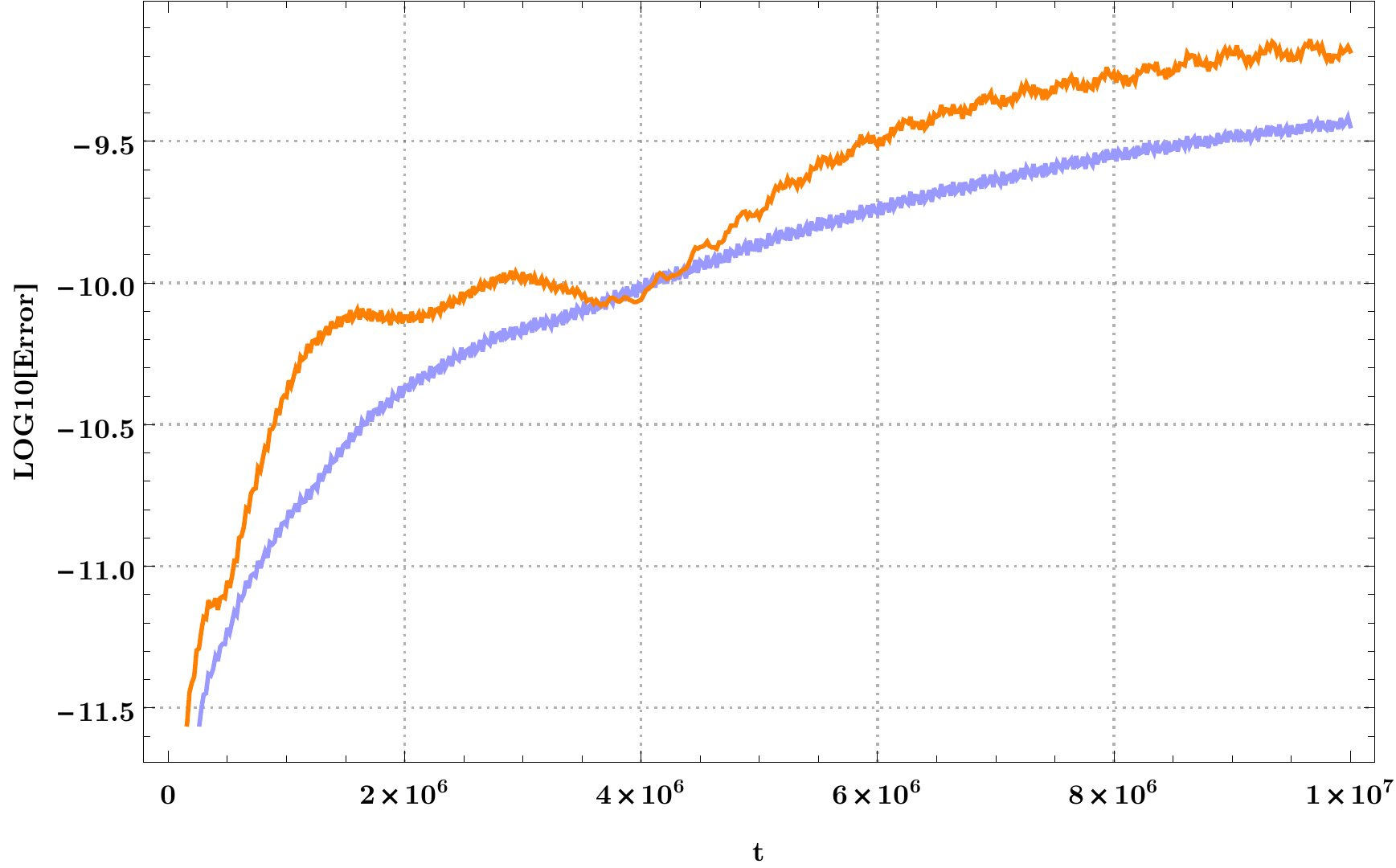}}
&
\subfloat[OSS: $P=1000$ perturbed initial values]
{\includegraphics[width=.4\textwidth]{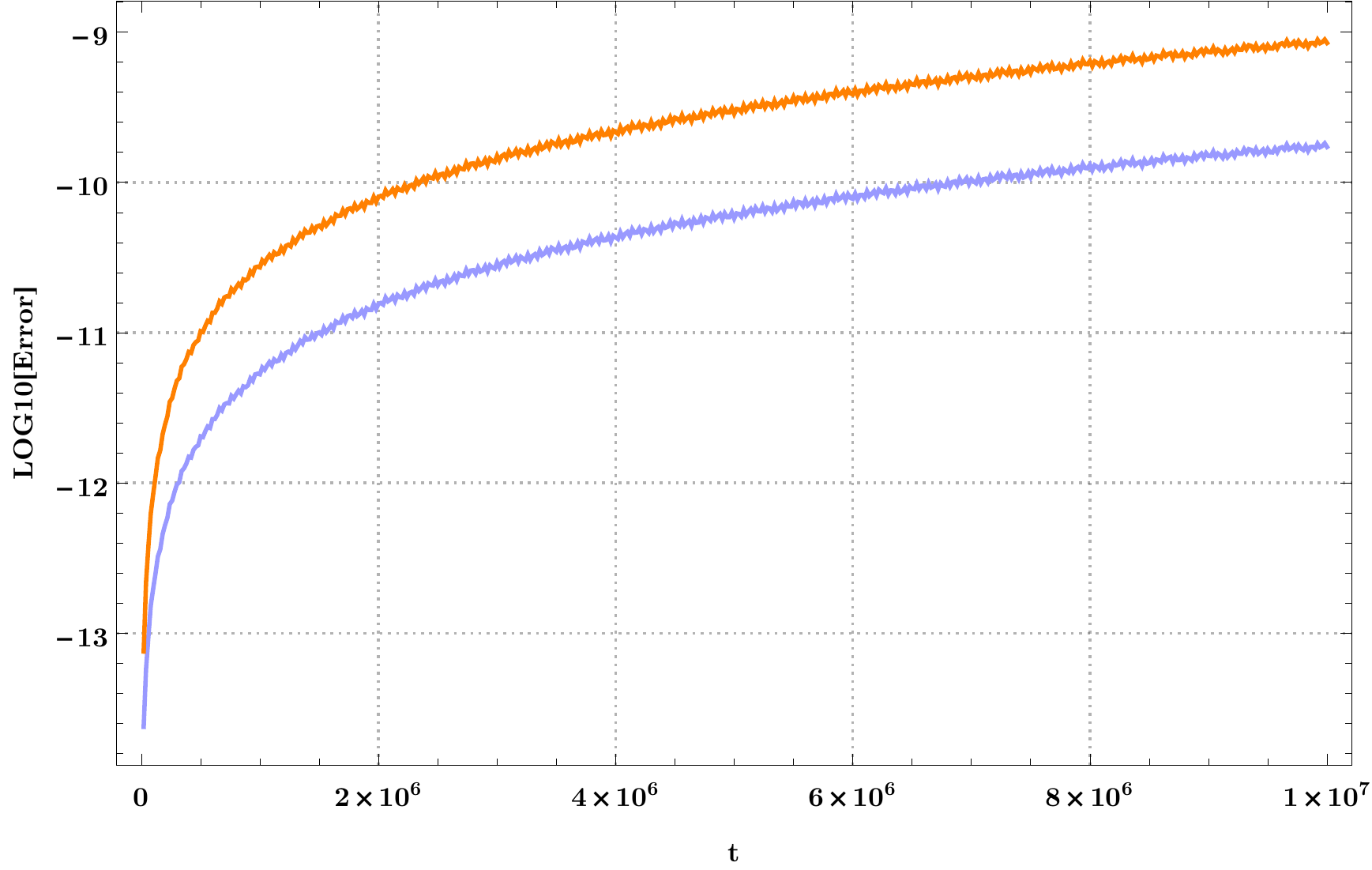}}
\end{tabular}
\caption{\small Left: estimation of the round-off error with the original unperturbed initial values. We compare the evolution of our error estimation (orange) with the evolution of the global error (blue). Right: evolution of the mean error in positions (blue) of the application of our DP algorithm to $P=1000$ perturbed initial value problems, together with the evolution of the mean of the estimated errors in positions (orange).}
\label{fig:plot5}
\end{figure}

In order to assess the quality of the error estimation technique proposed in Subsection~\ref{ss:estimation}, we represent, In Fig.~\ref{fig:plot5}, for each of the three considered initial value problems (with the original unperturbed initial values), the evolution of the global errors in position of our DP implementation, together with the evolution of the estimations produced by using our technique applied with $r=3$. In addition, we present for each of the three considered examples, the evolution of the mean error in positions of the application of our DP algorithm to $P=1000$ perturbed initial value problems, together with the evolution of the mean of the estimated errors in positions.
We  believe that the results indicate that the proposed round-off error estimation procedure is useful for the purpose of assessing the propagation of round-off errors.

\clearpage
\section{Concluding remarks}

Symplectic implicit Runge-Kutta schemes (such as RK collocation methods with Gaussian nodes) are very appropriate for the accurate numerical integration of general Hamiltonian systems. For non-stiff problems, implementations based on fixed-point iterations seem to be more efficient than those based on Newton method or some of its variants.

We propose an implementation that takes special care in reducing the propagation of round-off errors, and  includes the option of computing, in addition to the numerical solution, an estimation of the propagated round-off error. We claim that our implementation with fixed point iterations is near optimal, in the sense that the propagation of round-off errors is essentially no worse than the best possible implementation with fixed point iteration.  Our claim seems to be confirmed by our numerical experiments.

A key point in our implementation has been the introduction of a new stopping criterion for the fixed point iteration.  We believe that such a stopping criterion could be also useful in other contexts.

According to our numerical experiments, it seems that, in some cases,  some small linear drift of the mean energy error may be unavoidable for the fixed point implementations of IRK schemes.  Whenever avoiding  any drift of energy error becomes critical it might be preferable to use some Newton based iteration instead.

The C code of our implicit Runge-Kutta implementation with fixed point iterations can be downloaded from \href{<https://github.com/mikelehu/IRK-FixedPoint>}{IRK-FixedPoint}  Github software repository or go to the next url: \url{https://github.com/mikelehu/IRK-FixedPoint} .

\paragraph{Acknowledgements}

M. Anto\~nana, J. Makazaga, and A. Murua have been partially supported by projects MTM2013-46553-C3-2-P from Ministerio de Econom\'ia y Comercio, Spain, by project MTM2016-76329-R “IMAGEARTH”  from Spanish Ministry of Economy and Competitiveness
and as part of the Consolidated Research Group IT649-13 by the Basque Government.

\section*{Appendix A: Computation of coefficients for 12th order Gauss collocation method}

We next illustrate, by considering in detail the case of the 6th stage Gauss collocation method for the 64-bit IEEE double precision  floating point arithmetic,  how to determine appropriate machine number coefficients  $\mu_{i j}$, $1 \leq i,j \leq s$, that approximate the real numbers $a_{i j}/b_j$ of a given symplectic integration.

For all $i=1,\ldots,s$, $\mu_{i,i}=1/2$.  For $1 \leq j < i \leq s$, $ \mu_{i j}:= \fl(a_{i j}/b_{j})$,  which satisfy $1/2 < | \mu_{i j} | < 2$, which implies that $\mu_{j i}:= 1 -  \mu_{i j}$ is a machine number.  This results in machine number coefficients $\mu_{i j}$ that satisfy the symplecticity conditions (\ref{eq:sympl_cond_2}).
 
Given $h$, the coefficients $hb_i = h \times b_i$ are precomputed as follows: for $i=2,\ldots,s-1$, $hb_i := \fl(h \times b_i)$, and
\begin{equation*}
hb_1 := hb_s := (h - \sum_{i=2}^{s} hb_i)/2.
\end{equation*}

\end{document}